\documentclass[11pt]{article}

\usepackage[left=1.10in,right=1.10in,top=1in]{geometry}
\usepackage{latexsym, graphicx, epsfig, amsmath, amsfonts,amssymb,subfigure,color,float}
\usepackage{authblk}
\newenvironment{revision1}{\par}{\par}
\newenvironment{revision2}{\par}{\par}

\begin{document}

%\begin{frontmatter}

\title{An Efficient Solver for Cumulative Density Function-based Solutions of Uncertain Kinematic Wave Models}
%\tnoteref{mytitlenote}}


\author[1]{Ming Cheng}
\author[1]{Yi Qin}
\author[2]{Akil Narayan}
\author[3]{Xinghui Zhong}
\author[4]{Xueyu Zhu}
\author[1,5]{Peng Wang}%\corref{mycorrespondingauthor}}

%\ead{wang.peng@buaa.edu.cn, wang.peng@outlook.com}

\affil[1]{School of Mathematics and System Sciences, Beihang University, Beijing, China }
\affil[2]{Department of Mathematics, and Scientific Computing and Imaging (SCI) Institute, The University of Utah, Salt Lake City, UT}
\affil[3]{School of Mathematics, Zhejiang University, Hangzhou, China }
\affil[4]{Department of Mathematics, University of Iowa, Iowa, U.S.A.}
\affil[5]{Beijing Advanced Innovation Center for Big Data and Brain Computing, Beihang University, Beijing 100191, China}

\maketitle

\begin{abstract}
We develop a numerical framework to implement the cumulative density function (CDF) method for obtaining the probability distribution of the system state described by a kinematic wave model. The approach relies on Monte Carlo Simulations (MCS) of the fine-grained CDF equation of system state, as derived by the CDF method. This fine-grained CDF equation is solved via the method of characteristics. Each method of characteristics solution is far more computationally efficient than the direct solution of the kinematic wave model, and the MCS estimator of the CDF converges relatively quickly. We verify the accuracy and robustness of our procedure via comparison with direct MCS of a particular kinematic wave system, the Saint-Venant equation.
\end{abstract}

%\begin{keywords}
%CDF method, kinematic wave equation, uncertainty quantification
%\end{keywords}

%\end{frontmatter}

%\linenumbers
%%%%%%%%%%%%%%%%%%%%%%%%%%%%%%%%%%%%%%%%%%%%%%%%%%%%%%%%%%%%%%%%%%%%%%%%%%%%%%%%%%
%%%%%%%%%%%%%%%%%%%%%%%%%%%%%%%%%%%%%%%%% INTRODUCTION %%%%%%%%%%%%%%%%%%%%%%%%%%%%%%%%%%%%%%%%%
%%%%%%%%%%%%%%%%%%%%%%%%%%%%%%%%%%%%%%%%%%%%%%%%%%%%%%%%%%%%%%%%%%%%%%%%%%%%%%%%%%

%%%%%%%%%%%%%%%%%%%%%%%%%%%%%%%%%%%%%%%
\section{Introduction}\label{sec:introduction}
%%%%%%%%%%%%%%%%%%%%%%%%%%%%%%%%%%%%%%%

Kinematic wave models (KW) \cite{Lighthill_1955a_Kinematic,Lighthill_1955b_Kinematic} provide an important mathematical tool to describe environmental flows, such as overland flow and erosion~\cite{Singh_1996_kinematic,Singh_2001_kinematic}. Due to multiscale heterogeneity and insufficient site characterization, input parameters or fields into KW models often exhibit a high degree of uncertainty and are commonly modeled as random quantities. Such a probabilistic approach renders an otherwise deterministic kinematic wave equations as stochastic, and results in complete solutions that are in the form of probabilistic density functions (PDFs) or cumulative density functions (CDFs) of system states. This amounts to forward propagation of parametric uncertainty through the modeling process.

To quantify such predictive uncertainty associated with random parameters, one can employ various uncertainty quantification tools. Monte Carlo Simulations (MCS) is a popular approach and has been applied to obtain spatially-distributed probabilistic prediction from the stochastic kinematic wave model~\cite{Luc2002_Impact,morbidelli-2007-simplified}. It computes the mean and variances of system states by solving the deterministic governing equations multiple times using a number of realizations of the input random parameters or fields. Statistical moments can be then used as a forecast of the system average response and/or a measure of associated prediction error. However, if one is interested in the distribution tail of system states, a common concern in the risk assessment of many environmental applications~\cite{tartakovsky-2007-probabilistic}, the MCS framework becomes expensive due to the requisite large number of simulations for accurate estimation of these higher-order moments.

Another popular uncertainty quantification framework is the generalised polynomial chaos (gPC) expansion method~\cite{XiuK_SISC02} that aims to construct a surrogate relationship between system states and random parameters. Based on the seminal work on Hermite polynomial chaos~\cite{GhanemS91}, the gPC method is mathematically robust and is numerically easy to implement with stochastic Galerkin method or various stochastic collocation schemes~\cite{MathelinH_NASA03}, such as sparse grid method~\cite{XiuH_SISC05}, pseudospectral method~\cite{XiuS_CICP07}, compressive sensing method and various adaptive schemes~\cite{FrauenfelderST_CMAME05, WanK_JCP05, ZabarasG_JCP07}. However, the ``curse of dimensionality" is a major concern that incurs heavy computational cost when the uncertain random fields have weak two-point correlation functions as happens with, for example, white noise~\cite{Xiu_CICP09,Xiu10}.  

In recent years the so-called CDF method was proposed to address the parametric uncertainty of kinematic wave models~\cite{wang-2012-uncertainty}. Based on the concept of (one-point, one-time, Eulerian) fine-grained CDF, it is an extension of early works on the PDF method~\cite{Pope2000_turbulent, Shvidler2003_probability, tartakovsky-2011-pdf} and has found applications in two-phase flow in porous media~\cite{wang-2013-cdf, boso-2014-cumulative, Boso-2016-distribution} and others. Compared to direct simulations of the kinematic wave model, the CDF method offers three major advantages: (1) linearity of new fine-grained CDF equation can drastically reduce the computational cost; (2) one needs only the first statistical moment of solution to the fine-grained CDF equation and (3) its convergence rate is independent of the number of random parameters. Previous studies have since proposed semi-analytical solutions to special cases of the kinematic wave model~\cite{wang-2012-uncertainty}, but a comprehensive and efficient numerical scheme is yet developed for the fine-grained CDF equation. This is  the focus of our study.

In this paper, we present a numerical scheme to implement the CDF method for the kinematic wave equation. In Section~\ref{sec:formulation}, we provide a brief review of the kinematic wave model and then derive its corresponding fine-grained CDF equation using the CDF method. A comprehensive numerical scheme is then proposed and verified with a simple case in Section~\ref{sec:results}. We investigate the robustness and salient features of the numerical scheme in Section~\ref{sec:results} by comparing the CDF solutions of one kinematic wave system, the Saint-Venant equation, with those obtained from direct MCS of the original kinematic wave model. Key conclusions are drawn in Section~\ref{sec:summary}. 

%%%%%%%%%%%%%%%%%%%%%%%%%%%%%%%%%%%%%%%%%%%%%%%%%%%%%%%%%%%%%%%%%%%%%%%%%%%%%%%%%%
%%%%%%%%%%%%%%%%%%%%%%%%%%%%%%%%%%%%%%%%% INTRODUCTION %%%%%%%%%%%%%%%%%%%%%%%%%%%%%%%%%%%%%%%%%
%%%%%%%%%%%%%%%%%%%%%%%%%%%%%%%%%%%%%%%%%%%%%%%%%%%%%%%%%%%%%%%%%%%%%%%%%%%%%%%%%%

%%%%%%%%%%%%%%%%%%%%%%%%%%%%%%%%%%%%%%%
\section{Problem Formulation} \label{sec:formulation}
%%%%%%%%%%%%%%%%%%%%%%%%%%%%%%%%%%%%%%%

In this Section we formulate the general problem of stochastic kinematic wave (Section \ref{sec:thesetup}) and present the CDF formulation from previous work~\cite{wang-2012-uncertainty} (Section \ref{sec:oldwork}). Without specifying otherwise, our work is conducted on a complete probability space $(\Omega, \mathcal{F}, \mu)$, where $\Omega$ is the collection of events, $\mathcal{F}$ represents a $\sigma$-algebra on sets of $\Omega$, and $\mu$ denotes a probability measure on $\mathcal{F}$. We assume that all random variables have finite second moment. Vectors are represented by lowercase boldface letters, e.g., $\mathbf q$, and scalars by lowercase normal letters, .e.g., $k$. The superscript in $\mathbf q^{\mathrm T}$ represents the vector transpose operator. We reserve $z$ and all its variations ($\boldsymbol{z}, Z, \boldsymbol{Z}$) for random variables.

%%%%%%%%%%%%%%%%%%%%%%%%%%%%%%%%%%%%%%%%%%%%%%%%%%%%%
\subsection{Kinematic wave equation} \label{sec:thesetup}
%%%%%%%%%%%%%%%%%%%%%%%%%%%%%%%%%%%%%%%

Kinematic wave models are popular tools in the study of flow motions driven by gravity and pressure. These models are constructed by postulating functional relationships between the quantity of interest $k(\mathbf x,t)$ (quantity per unit distance) and its flux $\mathbf q(\mathbf x,t)$ (quantity passing a given point in unit time and distance), 
\begin{subequations} \label{model:KWT}
	\begin{equation} 
		\mathbf q= \mathbf q(k, \mathbf z), \label{eqn:con}
	\end{equation}
where $\boldsymbol{z}$ above is a random vector that encodes the uncertainty in this relationship induced by, e.g., uncertain environmental properties.

The kinematic wave model combines the conservation laws of mass and momentum and describes the wave phenomenon with the continuity equation alone:
	\begin{eqnarray} \label{KWT}
		\frac{\partial k}{\partial t} + \nabla \cdot \mathbf q = S, 
		\qquad \mathbf x \equiv(x_1, x_2, x_3)\in  \mathcal D, \; t>0,	
	\end{eqnarray}
subject to the following initial and boundary conditions:
	\begin{eqnarray}  
		\label{initial-k}
               & k_{\mathrm{in}}(\mathbf x) =  k(\mathbf x,t=0), \quad &\boldsymbol{x} \in \mathcal{D}, \\
		\label{boundary-k}
               & k_{\mathrm {bx}}(t) =  k(\mathbf x = \partial \mathcal{D} ,t), \quad &\boldsymbol{x} \in \partial\mathcal{D}.
	\end{eqnarray} 
\end{subequations}
Here, $S(\mathbf x,t)$ is a source/sink term, $\mathbf z\equiv(z_1, \dots, z_\mathrm{N})$ are $N$ random parameters in the constitutive relationship~\eqref{eqn:con}, and $\mathcal{D}$ is a spatial domain.

In general, the constitutive relationship~\eqref{eqn:con} and its parameterization $\mathbf z$ are often based on physical interpretations of the underlying process. For example, in overland flows one may use Darcy-Weisbach, Ch{\' e}zy, or Manning formulae to describe laminar, turbulent or transitional flow regime, respectively~\cite{Singh_2001_kinematic}.

In many environmental applications, source/sink, functional parameters $\mathbf z $, boundary and initial conditions may be subject to epistemic uncertainty. In the case of open-channel flow, $S(\mathbf x, t)$ represents the rainfall rate and/or inflow/outflow at the tributaries and/or runoff from the ambient terrain, all of which exhibit heterogeneity at various spatio-temporal scales. Meanwhile, past data analysis (\cite{Buhman-2002-Stochastic,Gates_1996a_Spatio,Gates_1996b_Spatio,Moramarco_2000_Kinematic} and the references therein) of the two functional parameters, e.g. the channel slope and surface roughness coefficient, suggest site-specific statistical distributions to capture their spatial variability. Although data acquisition continues to improve, ubiquitous data sparsity and measurement/interpretation errors render overland flow predictions inherently uncertain. This predictive uncertainty is routinely mentioned as one of the fundamental challenges in flood forecasting \cite{moore-2006-issues}.

In subsequent analysis, we model the parametric uncertainty as (un/correlated) random fields. In other words, a quantity of interest $k(\mathbf x, t, \omega)$ varies not only in the physical domain, $(\mathbf x,t)\in \mathcal{D} \times [0, \infty)$, but also with respect to the probability event $\omega \in \Omega$. 

%
%%%%%%%%%%%%%%%%%%%%%%%%%%%%%%%%%%%%%%%%%%%%%%%%%%%%%%%%%%
\subsection{CDF method} \label{sec:oldwork}
%%%%%%%%%%%%%%%%%%%%%%%%%%%%%%%%%%%%%%%
To solve the stochastic kinematic wave equation~\eqref{model:KWT} and obtain a statistical description of $k$ at any space-time point, one can employ the CDF method~\cite{wang-2012-uncertainty, wang-2013-cdf, boso-2014-cumulative}. In the deterministic setting, we introduce the concept of (one-point, one-time, Eulerian) fine-grained cumulative density function (CDF) of $k$:
%. Similar to the PDF method, as pioneered by Pope in the study turbulence~\cite{Pope2000_turbulent} and later in contaminant reactive transport in porous media \cite{Shvidler2003_probability, Tartakovsky_2003_PDF}, the 
%
\begin{equation}\label{def:Pi}
	\Pi(K; \mathbf x, t) = \mathcal{H}  [ K - k(\mathbf x, t)],
	\quad K \in \mathbb R^{+},
\end{equation}
where $\mathcal{H}(\cdot) $ is the Heaviside step function and $K$ is a deterministic value (outcome) that the random quantity $k$ can take at a space-time point $(\mathbf x,t)$. The Heaviside function is defined here as
\begin{align*}
  \mathcal{H}(y) = \left\{\begin{array}{rl} 0, & y < 0 \\ 1, & y \geq 0. \end{array}\right.
\end{align*}
In general, the relation \eqref{eqn:con} with random variables $\boldsymbol{z}$ induces randomness in the solution. In this case the function $\Pi$ is random, a function of $\omega$, but in the sequel we continue to refer to it as a ``CDF", consistent with previous literature. With $f_{\mathrm k}(K; \mathbf x, t)$ denoting the single-point probability density function (PDF) of $k(\mathbf x,t,\omega)$, at fixed $(\boldsymbol{x},t)$, the expectation of the $\Pi$ as defined in~\eqref{def:Pi} yields the single-point CDF $F_{\mathrm k}(K; \mathbf x, t)$,
\begin{equation} \label{eqn:piCDF}
 	\mathbb{E}\; \Pi(K; \mathbf x,t) := \int_{-\infty}^{\infty} \mathcal{H}  (K-k') f_ {\mathrm k}(k'; \mathbf x, t) {\mathrm d} k'  =  F_{\mathrm k}(K; \mathbf x,t).
\end{equation}
Here we emphasize again, that the ``CDF" $\Pi$ defined in \eqref{def:Pi} is an indicator function, and is random. The actual PDF and CDF of $k$ are $f_{\mathrm k}$ and $F_{\mathrm k}$, respectively, and are deterministic. Following earlier works~\cite{wang-2012-uncertainty}, we multiply the kinematic wave equation~\eqref{KWT} with $\partial \Pi / \partial K$ and obtain a linear fine-grained CDF equation (see Appendix~\ref{appendix:pi}):
\begin{subequations} \label{eqn:pi}
\begin{equation}
	 \frac{ \partial \Pi } {\partial t} +  \mathbf v(K; \mathbf{x}, t) \cdot \nabla_{\mathbf c} \Pi =0,
\end{equation}
in which the operator $\nabla_{\mathbf c}(\cdot)$ and advection velocity $\mathbf v$ are:
\begin{eqnarray}
		&& \qquad \qquad 
	\nabla_{\mathbf {c}} = \left( 
	\frac{\partial}{\partial x_1}, \, \frac{\partial}{\partial x_2}, \, \frac{\partial}{\partial x_3}, \, \frac{\partial}{\partial K} \right)^\mathrm{T},\\
	 && \mathbf {v} (K; \mathbf{x}, t) =  \left(
%	\begin{array}{c}
 \frac{\partial q_1}{\partial K} , \;
 \frac{\partial q_2}{\partial K} , \;
 \frac{\partial q_3}{\partial K}, \;
\sum_{i, j=1}^{3, \mathrm N} \frac{\partial q_i}{\partial z_j} \, \frac{\partial z_j}{\partial x_i} + S 
%\vdots \\
 %	\end{array}
 	\right) ^\mathrm{T},
\end{eqnarray}
\end{subequations}
whose initial and boundary conditions are derived from the physical space relations~\eqref{initial-k}\&\eqref{boundary-k}:
\begin{subequations}\label{bc-PI}
\begin{eqnarray} 
 & \Pi_\mathrm{in} =  \Pi(K; \mathbf x, 0)  = \mathcal{H}  \left[ K- k_\mathrm{in}(\mathbf x) \right], & \boldsymbol{x} \in \mathcal{D} \label{eqn:bc:pit} \\ 
&           \Pi_\mathrm{bx}  = \Pi(K; \boldsymbol{x},t)= \mathcal{H}  \left[ K - k_\mathrm{bx}(t) \right], & \boldsymbol{x} \in \partial \mathcal{D}\label{eqn:bc:pix} 
	   %\quad \forall \mathbf x \in \partial \mathcal D,
\end{eqnarray}
For the additional dimension $K$, a boundary condition~$\Pi_\mathrm{bk}$ must be prescribed for a unique solution of $\Pi(K; \mathbf x, t)$. It can be determined often intuitively from specific conditions of the underlying physical process, for example, in open-channel flow, the water height $k$ may be assumed always greater than zero:
\begin{eqnarray}
	 \Pi _\mathrm{bk} =\;   \Pi(0; \mathbf x,t)  =0.  \label{eqn:bc:piK}
\end{eqnarray}
\end{subequations}

The CDF formulation offers four major advantages. First, comparing to the nonlinear kinematic wave equation~\eqref{KWT}, the linear fine-grained CDF equation~\eqref{eqn:pi} is easier to solve, albeit in a higher spatial dimension. Second, by invoking~\eqref{eqn:piCDF}, one only needs to compute the ensemble average of $\Pi$ to obtain the full statistical distribution, bypassing the need to compute high-order moments. Third, the CDF formulation does not impose any prior assumption on the number of random parameters or on their correlation structures. Lastly,  the fine-grained CDF boundary condition \eqref{eqn:bc:piK} is formulated in a straightforward and unambiguous way, whereas its counterpart in the PDF methods~\cite{Pope2000_turbulent,tartakovsky-2011-pdf}, e.g. the probability density of a specific condition, may be hard to formulate or is generally unknown.

%%%%%%%%%%%%%%%%%%%%%%%%%%%%%%%%%%%%%%%%%%%%%%%%%%%%%%%%
%%%%%%%%%%%%%%%%%%%%%%      METHOD       %%%%%%%%%%%%%%%%%%%
%%%%%%%%%%%%%%%%%%%%%%%%%%%%%%%%%%%%%%%%%%%%%%%%%%%%%%%%

%%%%%%%%%%%%%%%%%%%%%%%%%%%%%%%%%%%%%%%
\section{Numerical scheme for the CDF method} \label{sec:method}
%%%%%%%%%%%%%%%%%%%%%%%%%%%%%%%%%%%%%%%

In this section, we present an efficient numerical scheme (Section~\ref{subsec:numerical}) for the CDF formulation to obtain a complete (single space-time point) probabilistic description of $k(\mathbf x, t)$ via its cumulative density function $F_{\mathrm k}(K; \, \mathbf x, t)$. %This is accomplished through the Monte Carlo estimator \eqref{eq:F-mcs} via repeated solves of \eqref{eqn:pi}. 
We then consider a few examples, namely, one-dimensional flow, three-dimensional flow, coupled system of one-dimensional flow and Burgers' equation, to examine the accuracy and illustrate some salient properties of the numerical scheme (Section~\ref{subsec:example}).

%%%%%%%%%%%%%%%%%%%%%%%%%%%%%%%%%%%%%%%%%%%%%%%%%%%%%%%%%
\subsection{Numerical scheme} \label{subsec:numerical}

We propose an efficient numerical framework that focuses on the computation of the fine-grained CDF equation~\eqref{eqn:pi}. Due to its linearity, we can use the well-known method of characteristics to solve the $\Pi$ equation. Specifically, we obtain the family of characteristics,
\begin{equation} \label{eqn:cl-v}
	\frac{\mathrm d \mathbf c}{\mathrm d t} = \mathbf v,
	\quad \mathbf c= (x_1, x_2, x_3, K)^{\mathrm T} ,
	\qquad \mathbf c(t=0) = \mathbf c_0.
\end{equation}
along which the fine-grained CDF equation~\eqref{eqn:pi} reads:
\begin{equation} \label{eqn:cl-k}
	\frac{\mathrm d \Pi}{\mathrm d t} = 0.
\end{equation}
The original problem~\eqref{eqn:pi} is therefore recast to a set of ordinary differential equations~\eqref{eqn:cl-v} and \eqref{eqn:cl-k}, which can be then solved by standard numerical ODE solvers.  We employ a third-order Total Variation Diminishing (TVD) Runge-Kutta scheme~\cite{Gottlieb-1998-Total}.

\begin{revision1}
By invoking \eqref{eqn:piCDF}, we can create an ensemble of $M$ realization of $\boldsymbol{z}$ in \eqref{eqn:con}, compute the associated ensemble of solutions $\Pi$ by solving \eqref{eqn:pi}, and finally compute $F_{\mathrm k,M}(K;\boldsymbol{x},t)$, which is an approximation to integral $F_\mathrm{k}(K; \mathbf x, t)$ from the ensemble of its $\Pi$ solutions. 

\begin{eqnarray}\label{eq:F-mcs}
	F_{\mathrm k}(K;\boldsymbol{x},t) &=& \int_{-\infty}^{\infty} \Pi \, \mathrm d F_ {\boldsymbol z}(\boldsymbol z'; \mathbf x, t) \nonumber \\
	&\approx& \sum_{i=0}^M \Pi_i(\boldsymbol z_i, K, \boldsymbol{x}, t) \left[ F_{\boldsymbol z}( \boldsymbol z_{i+1}') - F_{\boldsymbol z}( \boldsymbol z_i')\right] 
	=F_{\mathrm k,M}(K;\boldsymbol{x},t), 
%	&& \begin{cases}
%	%\left\{
%	i=0, F_{\boldsymbol z} (\boldsymbol z_{0}') = 0\\
%	i=M, F_{\boldsymbol z} (\boldsymbol z_{M+1}') = 1
%	%\right. 
%	\end{cases}
%	\nonumber 
\end{eqnarray}
where $\boldsymbol{z}_i$ represents the $i$-th realization of $\boldsymbol{z}$, $\Pi_i$ is the corresponding fine-grained CDF solution, $F_ {\boldsymbol z}(\boldsymbol z_i')$ is the joint probability of random inputs for that realization while $ F_{\boldsymbol z} (\boldsymbol z_{0}') = 0$ and $F_{\boldsymbol z} (\boldsymbol z_{M+1}') = 1$, respectively. For independent random variables $\boldsymbol z$, the joint probability $F_ {\boldsymbol z, i}$ would be a tensor product of each random inputs and various collocation techniques can be implemented. For correlated random inputs, one may employ the latest work~\cite{Akil-2018-high-dimension} to compute the high dimensional integral with designed quadrature. We hope to address this subject in future works. %It is noted here that direct ensemble average of the fine-grained CDF solutions can also be adopted to relax any prior assumptions of the random parameters, which is in effect a MCS average. 
\end{revision1}

The overall procedure can be summarized as follows:
\begin{enumerate}
        \item{\it Generate $M$ realizations of the random parameters $\boldsymbol{z}$.}
	
	\item{\it Solve the fine-grained CDF equation~\eqref{eqn:piCDF} $M$ times via the method of characteristics, one for each realization.}
	
        \item{\it Compute the ensemble average of the $\Pi$ solutions using \eqref{eq:F-mcs}.}
\end{enumerate}
In essence, our numerical framework is based on two key features of the CDF formulation: linearity of the fine-grained CDF equation~\eqref{eqn:pi} and the ensemble average of its solution is the CDF of system state~\eqref{eqn:piCDF}. To implement the method of characteristics, one can take advantage of the large body of literature on the numerical methods for ordinary differential equations.

Alternatively, one may employ the Reynolds decomposition: $\mathcal A= \langle \mathcal A \rangle+ \mathcal A'$, to represent the random parameters as the sum of their ensemble means  $\langle \mathcal A \rangle$ and zero-mean fluctuations about the mean $\mathcal A'$. Then, by taking the ensemble average of the linear fine-grained equation~\eqref{eqn:pi}, a deterministic equation ensues:
\begin{eqnarray} \label{eq:advec}
 	\frac{\partial F_{\mathrm k}}{\partial t} + \mathbf v_\mathrm{eff} \cdot \nabla_{\mathbf c} F_\mathrm{k} =  \nabla_{\mathbf c} \cdot (\mathbf D \nabla_{\mathbf c} F_\mathrm{k}),
	% \quad \nabla_{\mathbf xK} = \left( \frac{\partial}{\partial x_1}, \frac{\partial}{\partial x_2}, \frac{\partial}{\partial x_3}, \frac{\partial}{\partial K}\right) ^\mathrm{T}
\end{eqnarray}
where $\mathbf v_\mathrm{eff}$ and $\mathbf D$ are the effective velocity and the eddy-diffusivity tensor, respectively. Similar to standard PDF methods~\cite{tartakovsky-2009-probability,tartakovsky-2011-pdf, wang-2013-probability, venturi-2013-exact}, such numerical approach~\eqref{eq:advec} requires a closure approximation, such as the large-eddy-diffusivity closure, and its solution is asymptotically exact when $F_\mathrm{k}$ varies slowly with $\mathbf x$, $K$ and $t$, relative to $\mathbf v$~\cite{kraichnan-1987-eddy, Barajas-2018-Probability}.

%%%%%%%%%%%%%%%%%%%%%%%%%%%%%%%%%%%%%%%%%%%%%%%%%%%%%%%%%
%%%%%%%%%%%%%%%%%%%%%%%%%%%%%%%%%%%%%%%%%%%%%%%%%%%%%%%%%
\subsection{Numerical examples}\label{subsec:example}
In this section, we present several numerical examples to illustrate the applicability and efficiency of our method.
%%%%%%%%%%%%%%%%%%%%%%%%%%%%%%%%%%%%%%%%%%%%%%%%%%%%%%%%%
\subsubsection{One-dimensional flow}
Let us consider the following kinematic wave model:
\begin{eqnarray}\label{eqn:test-kw}
	 \frac{\partial k}{\partial t} + \frac{\partial q}{\partial x} = S , \quad
	q=k^{\frac{1}{2}},
	\qquad x \in  [0,2], \, t \in [0, + \infty)
\end{eqnarray}
%\annote{What is the spatial domain here? $\mathcal{D} = [0,1]$? the domain of space x is $[0,2]$, and the domain of time is $[0,0.1]$ }
%
whose corresponding fine-grained CDF equation is:
\begin{equation} \label{eqn:test-pi}
	\frac{\partial \Pi}{\partial t} + \frac{1}{2\sqrt{K}} \frac{\partial \Pi}{\partial x} + S \frac{\partial \Pi}{\partial K}= 0.
\end{equation}

\paragraph{Deterministic case}
We first examine  the accuracy of our numerical scheme  by considering deterministic source term and initial and boundary conditions:
\begin{equation}
	\begin{split}
	&  S(x, t)= 2\pi \left[ \sin(\pi(x+t)) + 1.1\right] \, \cos\left( \pi (x+t)\right)  + \pi \left[ \cos\left( \pi(x+t)\right) \right], \\
	& \qquad \qquad k_{\mathrm {bx}} = \left[ \sin{\pi t} + 1.1 \right]^2,
	 \qquad \qquad k_{\mathrm{in}}= \left[ \sin{\pi x} + 1.1 \right]^2, \\
	\end{split}
\end{equation}
which yields an explicit closed-form solution:
\begin{equation} \label{sol:test1}
  k_{\mathrm{exact}}(x,t) = \left[ \sin\left( \pi (x+t)\right) + 1.1\right]^2.
\end{equation}

Now we implement our numerical framework and compute the fine-grained CDF equation~\eqref{eqn:test-pi} with a third-order Runge-Kutta scheme (Appendix~\ref{appendix:RK}). We compared the numerical results with the analytical solution~\eqref{sol:test1}. The error  is measured
\begin{align}\label{eq:eps-def}
  \epsilon := \frac{1}{N} \sum_{i=1}^N \left[  k_{\mathrm{exact}}(x_i, t) - k(x_i, t) \right]^2,
\end{align}
where $k(x_i, t)$ is the solution computed using the method of characteristics, and $x_i$ is an $N$-point equidistant grid on $[0,2]$.

Table~\ref{tab1:test-1-pi} shows the comparison results for the solution at time $t=0.1$ and one can see that, with a diminishing time-step (from $0.1$ to $0.0125$), the error (in $l_2$-norm) $\epsilon$ between our numerical result and the exact solution diminishes to $1.71 \times 10^{-5}$. We observe third-order convergence in time for our approach, which matches expectations for this third-order Runge-Kutta scheme. % In the same Table~\ref{tab1:test-1-pi}, we also include a comparison between the exact solution and a direct numerical result of the kinematic wave equation~\eqref{eqn:test-kw} using a fifth-order WENO-Roe scheme. This is done so as a validation of the direct MCS benchmark result  in order to validate its , which would be later used as a MCS scheme
\begin{table}[h]
\centering
\begin{tabular}{lcc}
\hline
  $\Delta t$  & $\epsilon $ & Convergence Rate  \\
\hline
$0.1$          &$8.95 \times 10^{-3}$     & --- \\
$0.05$         &$1.11 \times 10^{-3}$     & $3.02$     \\
$0.025$        &$1.37 \times 10^{-4}$     & $3.01$   \\
$0.0125$        &$1.71 \times 10^{-5}$     & $3.00$   \\
\hline
\end{tabular}
\caption{$l_2$ error between the numerical result and the exact solution~\eqref{sol:test1} at time $t=0.1$.}
\label{tab1:test-1-pi}
\end{table}
%

%%%%%%%%%%%%%%%%%%%%%%%%%%%%%%%%%%%%%%%%%%%%%%%%%%%%%%%%%
\paragraph{Stochastic case}
We now consider a random source but maintain deterministic initial and boundary conditions, 
\begin{equation}
	\begin{split}
	&  S(x, t)= 2 z^2 \pi \left[ \sin\pi(x+t) + 5 \right] \, \cos\pi (x+t)  + z \pi \, \cos \pi(x+t) , \\
	& \qquad \qquad k_{\mathrm {bx}} = \left[ \sin{\pi t} + 5 \right]^2,
	 \qquad \qquad k_{\mathrm{in}}= \left[ \sin{\pi x} + 5 \right]^2, \\
	\end{split}
\end{equation}
where $z$ is a random variable with lognormal distribution:
\begin{eqnarray}
	f_{\mathrm z} (Z; \mu, \sigma^2) = \frac{ 1 }{Z \sigma \sqrt{2 \pi}} \mathrm e^{- \left( \ln Z - \mu \right)^2/ 2 \sigma^2},
\end{eqnarray}
whose mean and variance are: $\mu=0$, $\sigma^2=0.1$.

The exact solution of the kinematic wave~\eqref{eqn:test-kw} is:
\begin{equation} \label{sol:test2}
  k(x,t,z) = \left[ z \sin\pi (x+t) + 5\right]^2,
\end{equation}
whose exact CDF can be found as:
\begin{eqnarray} \label{sol:test2-CDF}
	F_\mathrm{k}(K; x, t) &=& \frac{ f_{\mathrm z} \left[  (\sqrt{K}-5)/ \sin \pi (x+t)\right] + f_{\mathrm z} \left[  (-\sqrt{K}-5)/ \sin \pi (x+t)\right] }{2\sqrt{K}  \sin \pi (x+t)  }.
\end{eqnarray}

\begin{figure}[H]
\centering
\includegraphics[height=2.1in]{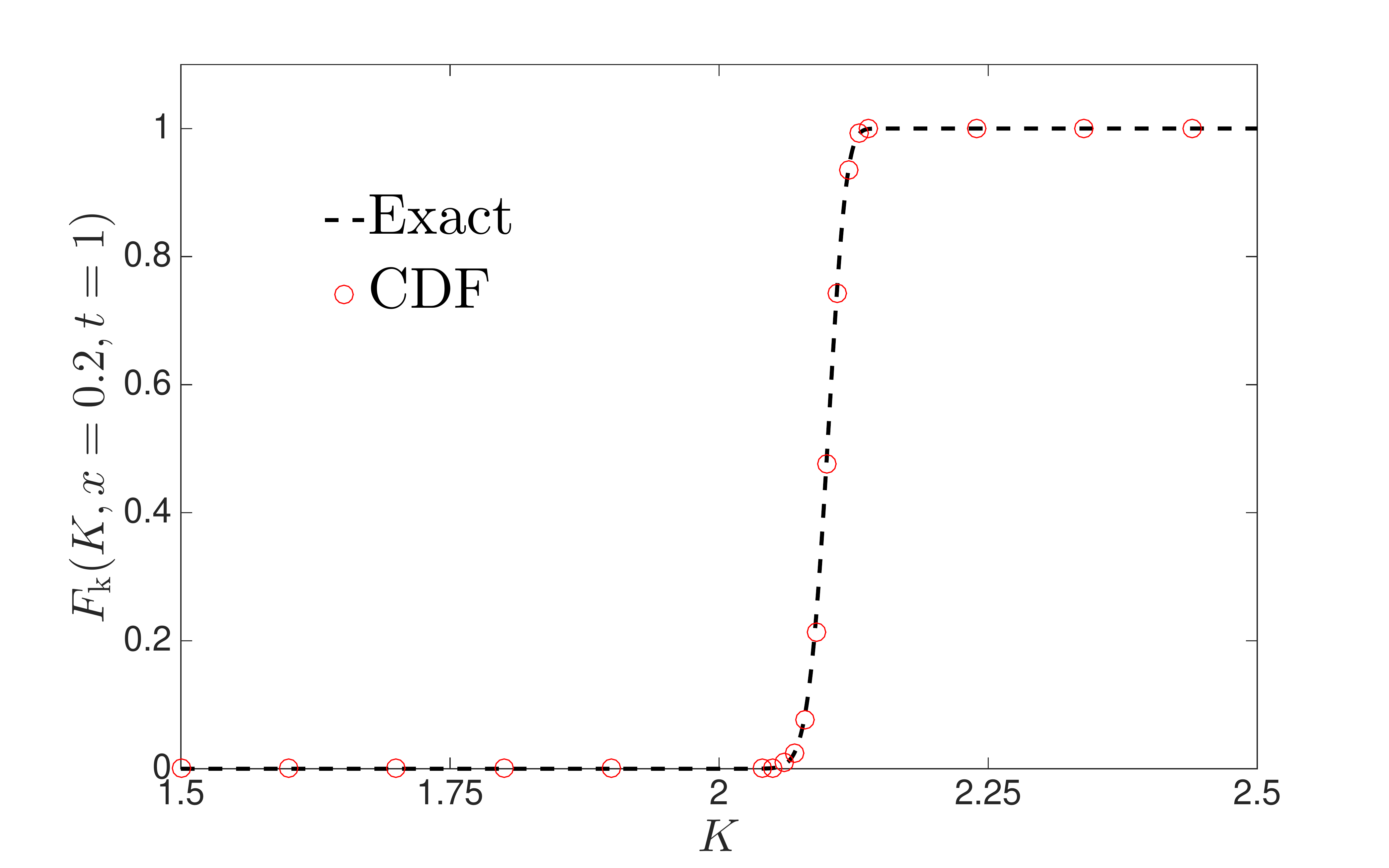}
\includegraphics[height=2.1in]{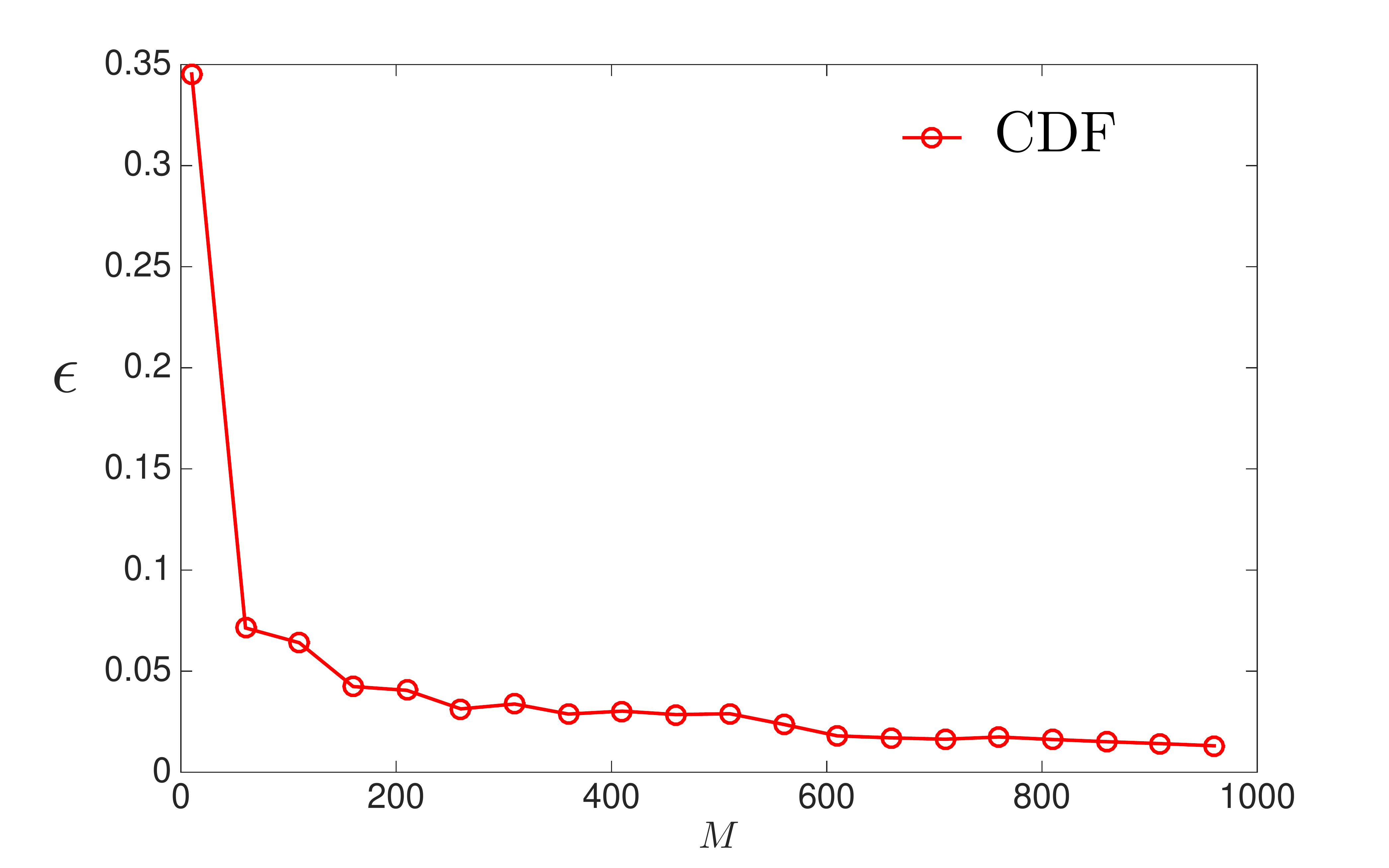}
  \caption{(a) Cumulative density function $F_{\mathrm k}(K; x=0.2, \, t=1)$ from the proposed CDF method $(CDF)$ and the exact $(exact)$ solution~\eqref{sol:test2-CDF}; (b) their relative error $\epsilon$ \eqref{error:1D-random} for different number of realizations.}
\label{fig:test2-solution}
\end{figure}

Using the proposed numerical scheme, we solve the stochastic fine-grained CDF equation~\eqref{eqn:test-pi} for a number of realizations of $z$. For each realization, it is found that the errors between exact and numerical solution are in the same order, e.g. $\sim O (10^{-5})$, as shown in the deterministic case. We then compute the ensemble average $F_{\mathrm{k},M}$ based on ~\eqref{eq:F-mcs} . In Fig.~\ref{fig:test2-solution} (a), the numerical CDF solution at a single space-time point, $F_{\mathrm k,M}(K; x=0.2, \, t=1)$, shows a good agreement with the analytical solution~\eqref{sol:test2-CDF}. For a closer examination in Fig.~\ref{fig:test2-solution} (b), we see that as the realizations number increases, the error ($\epsilon$) rapidly converges and then saturates after $100$ realizations. The error $\epsilon$ here for $F_{\mathrm k}$ is defined as:

\begin{align} \label{error:1D-random}
  \epsilon(x,t) := \frac{\sum_{i=1}^N \left[ F_{{\mathrm k},M}(K_i;x,t)- F_{\mathrm{k}}(K;x,t) \right]^2}{\sum_{i=1}^N F^2_{\mathrm{k}}(K;x,t)},
\end{align}
where $K_i$ is an $N$-point equidistant grid in $K$.

%%%%%%%%%%%%%%%%%%%%%%%%%%%%%%%%%%%%%%%%%%%%%%%%%%%%%%%%%
%%%%%%%%%%%%%%%%%%%%%%%%%%%%%%%%%%%%%%%%%%%%%%%%%%%%%%%%%\
\begin{revision1}
\subsubsection{ Three-dimensional flow}
We now consider a three-dimensional problem:
\begin{eqnarray} \label{eq:3d}
	&& \frac{\partial k}{\partial t}+\frac{\partial k}{\partial x_1}+\frac{\partial k}{\partial x_2}+\frac{\partial k}{\partial x_3} =  \frac{\pi z}{3} \left[ t\cos(\pi t x_1) + t \cos(\pi t x_2) \right. \nonumber \\
	&& \qquad \left.  + t \cos(\pi t x_3)  + x_1\cos(\pi t x_1) + x_2\cos(\pi t x_2) + x_3\cos(\pi t x_3) \right],
\end{eqnarray}
where $z$ is a random variable with lognormal distribution: $\ln z \sim \mathcal N (0,0.1^2)$.
One can find an exact solution to the three-dimensional flow as follows:
\begin{equation} \label{sol:3d}
	k=\frac{z}{3}\left[ \sin(\pi tx_1)+\sin(\pi tx_2)+\sin(\pi tx_3) \right].
\end{equation}
In Fig.~\ref{fig:3dsolution}, we plot the $l_2$-norm error~\eqref{error:1D-random} of the CDF solutions $F_{\mathrm k}(K; x_1=1.3, x_2=1.3, x_3=1.3, t=1)$ between the exact solution~\eqref{sol:3d} and the one obtained from the proposed CDF scheme, and that from the Monte Carlo Simulations (MCS) of the original system~\eqref{eq:3d}, respectively, at different realization number $M$. We found that the CDF method exhibits excellent accuracy with a fast convergence rate.
\begin{figure}[H]
\centering
\includegraphics[height=2.1in]{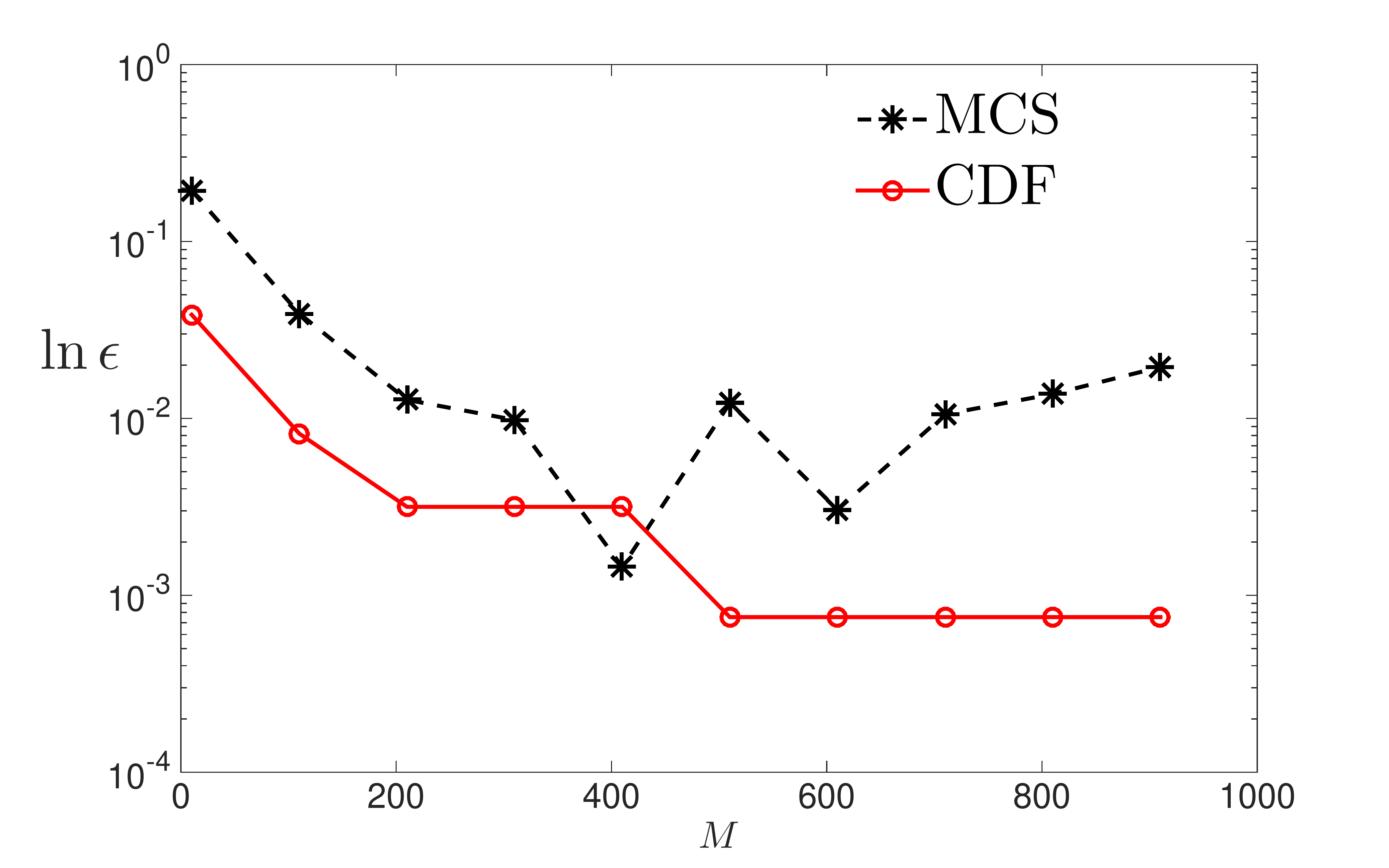}
  \caption{$l_2$-norm error~\eqref{error:1D-random} of the CDF solutions $F_{\mathrm k}(K; x_1=1.3, x_2=1.3, x_3=1.3 \, t=1)$ between the exact solution~\eqref{sol:3d} and the one obtained from the proposed CDF scheme, and that from the Monte Carlo Simulations (MCS) of the original system~\eqref{eq:3d}, respectively, at different realization number $M$.}
\label{fig:3dsolution}
\end{figure}

%%%%%%%%%%%%%%%%%%%%%%%%%%%%%%%%%%%%%%%%%%%%%%%%%%%%%%%%%
%%%%%%%%%%%%%%%%%%%%%%%%%%%%%%%%%%%%%%%%%%%%%%%%%%%%%%%%
\subsubsection{Coupled system of one-dimensional flow}
We now examine the following system of one-dimensional flow:
\begin{align} \label{eq:couple}
  \frac{\partial k_1}{\partial t}+\frac{\partial k_2}{\partial x} &=0,  \\
  \frac{\partial k_2}{\partial t}+\frac{\partial k_1}{\partial x} &=0, 
\end{align}
subject to the random initial conditions:
\begin{align} \label{sol:Eqs}
   k_1(x,0)&=z \sin(\pi tx),  \\
   k_2(x,0)&=z \cos(\pi tx),
\end{align}
where $z$ is a random variable with lognormal distribution: $\ln z \sim \mathcal N (0,0.1^2)$. 

Its exact solution can be found as
\begin{eqnarray}\label{sol:couple} 
	k_1 &=&\frac{1}{2} \left[ \sin\pi(x-t)+  \sin\pi(x+t)+  \cos\pi(x-t)- \cos\pi(x+t) \right], \\
k_2 &=&\frac{1}{2} \left[  \sin\pi(x-t)- \sin\pi(x+t)+ \cos\pi(x-t)+\cos\pi(x+t) \right].
\end{eqnarray}

We introduce two fine-grained CDFs: 
\begin{eqnarray}
	\Pi_{\mathrm v_1} = \mathcal H(v_1 - v_1'), 
	\quad
	\Pi_{\mathrm v_2} = \mathcal H(v_2-v_2'),
\end{eqnarray}
where the variables: $v_1=(k_1+k_2)/2$ and $v_2=(k_1-k_2)/2$ can help decouple the original system. Now one can employ the proposed CDF scheme to obtain those two marginal fine-grained CDFs for each realization and then recover the fine-grained CDFs of original system states: 
\begin{eqnarray}
	\Pi_{\mathrm k_1} &=& \mathcal H(k_1 - k_1') = \mathcal H\left[ (v_1 + v_2) \, - \, (v_1'+v_2')  \right], 
	\\
	\Pi_{\mathrm k_2} &=& \mathcal H(k_2 - k_2')= \mathcal H \left[ (v_1 - v_2) \, - \, (v_1' - v_2')  \right].
\end{eqnarray}
With a number of realizations, the marginal cumulative density functions $F_{\mathrm k_1}$ and $F_{\mathrm k_2}$ can be obtained with the approximation formula~\eqref{eq:F-mcs}. The $l_2$-norm errors~\eqref{error:1D-random} between the exact solutions~\eqref{sol:couple} and those from the CDF method, and MCS of the original coupled system, respectively, are plotted in Figure.~\ref{fig:couplesolution} for $F_{\mathrm k_1}(k_1'; x=0.3 , t=1)$ and $F_{\mathrm k_2}(k_2'; x=0.3 , t=1)$ at different realizations number $M$. It is clear that the CDF scheme provides better accuracy for the same number of realizations and achieves faster convergence rate. 

\begin{figure}[H]
\centering
\includegraphics[height=2.1in]{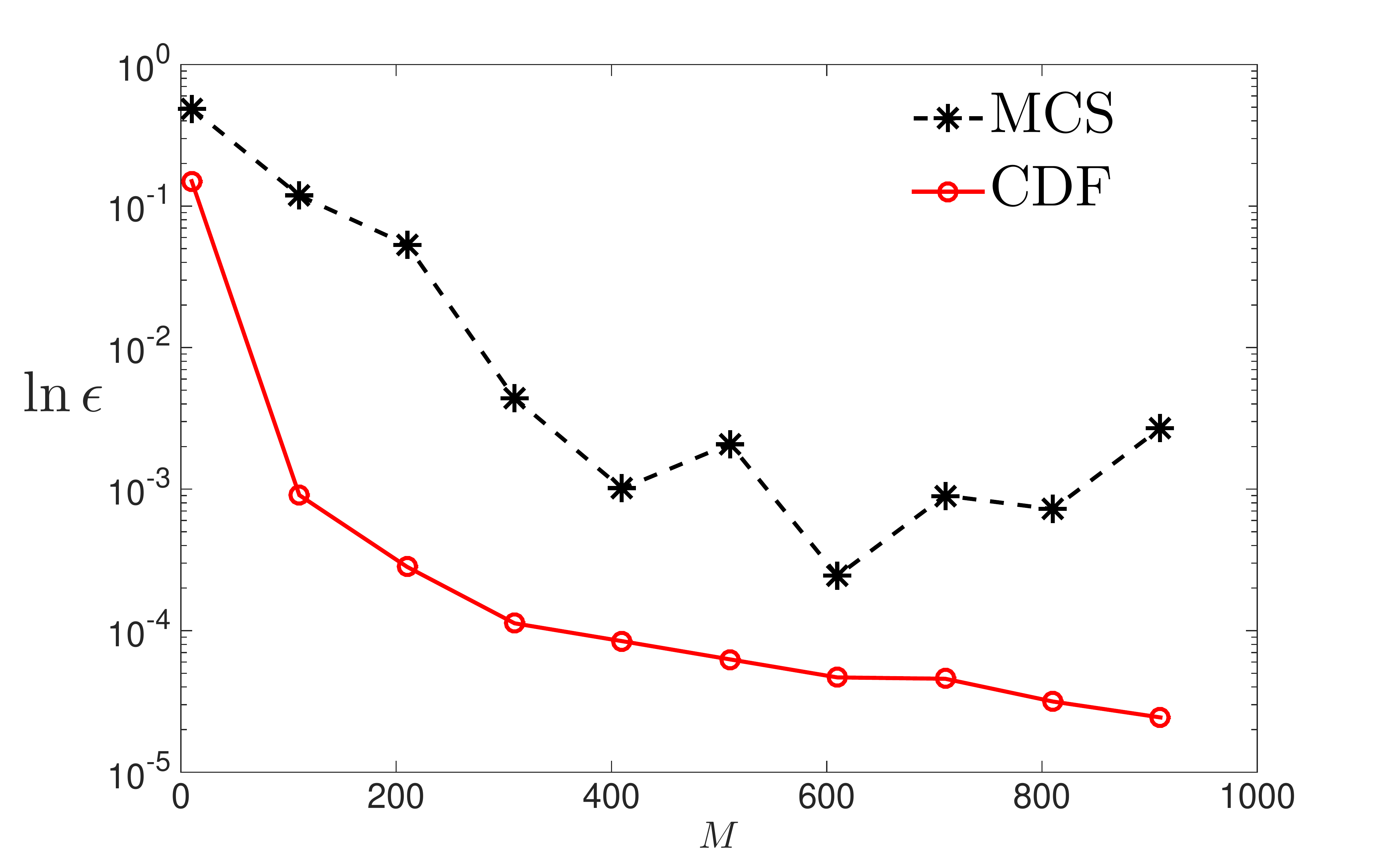}
\includegraphics[height=2.1in]{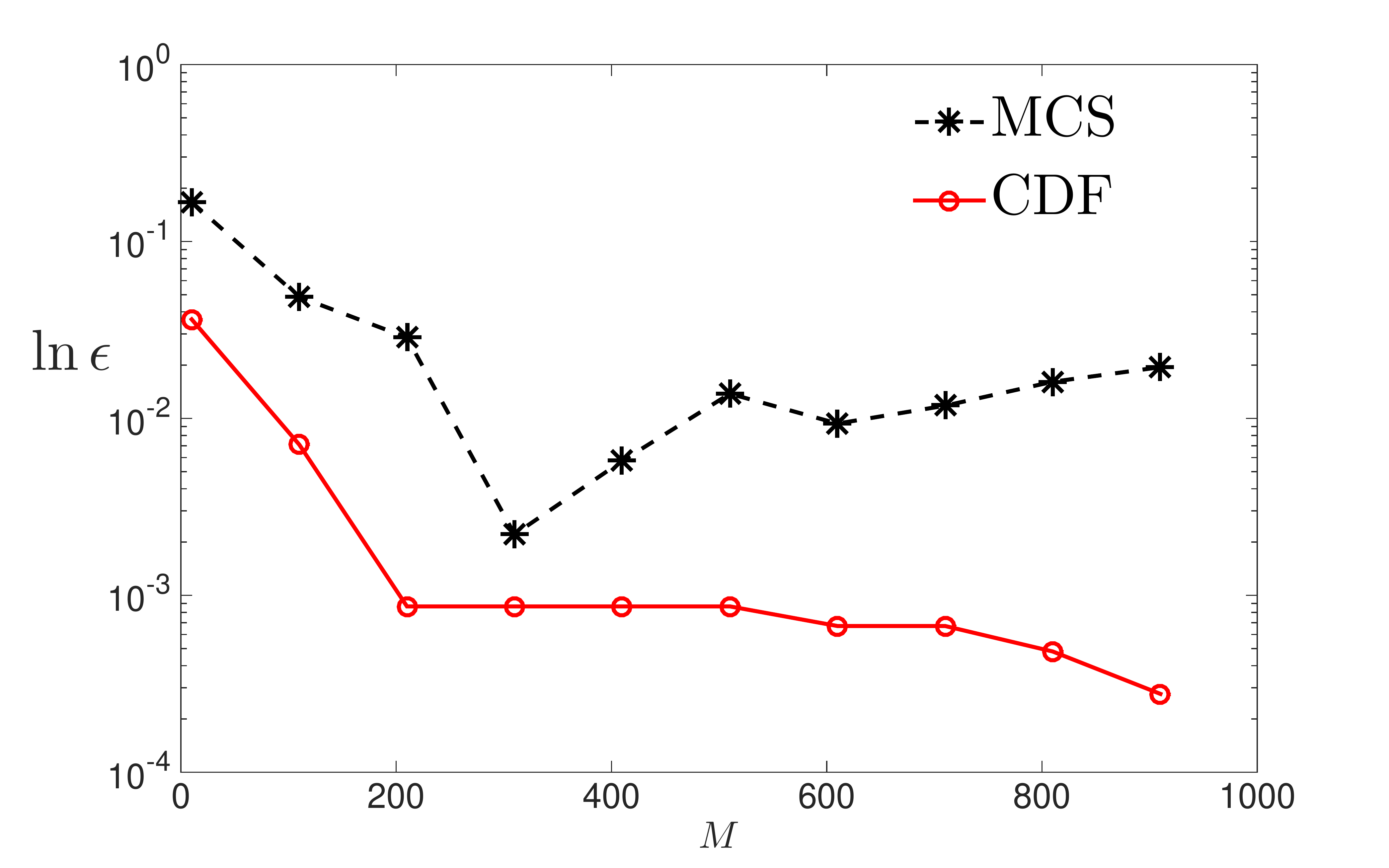}
  \caption{The $l_2$-norm errors~\eqref{error:1D-random} of a) $F_{\mathrm k_1}(k_1'; x=0.3 , t=1)$ and b) $F_{\mathrm k_2}(k_2'; x=0.3 , t=1)$, between the exact solutions~\eqref{sol:couple} and those from the CDF method, and MCS of the original coupled system, respectively.}
\label{fig:couplesolution}
\end{figure}
\end{revision1}
%%%%%%%%%%%%%%%%%%%%%%%%%%%%%%%%%%%%%%%%%%%%%%%%%%%%%%%%%%
%%%%%%%%%%%%%%%%%%%%%%%%%%%%%%%%%%%%%%%%%%%%%%%%%%%%%%%%%%
\begin{revision2}
\subsubsection{Burgers' equation}
Finally, let us consider a nonlinear example, the Burgers' equation:
\begin{equation} \label{eq:BG}
	\frac{\partial k}{\partial t}+\frac{\partial k^{2}}{\partial x}=0,
\end{equation}
subject to a stochastic initial condition with lognormal random variable $z$: $\ln z \sim \mathcal N (0,0.1^2)$.
\begin{equation} \label{sol:BG}
	k(x,0)=z\sin(\pi tx).
\end{equation}

The initial condition of the Burgers' equation may lead to shocks at later time. To address such issue, we follow earlier works~\cite{wang-2013-cdf,alawadhi-2018-method} and employ the Rankine-Hugoniot condition~\cite{Smoller-1983-Shock} to determine the shock location $x_{\mathrm s}$ at each realization. Here the fine-grained CDF solution would be divided by two parts: behind $x_{\mathrm s}$, the $\Pi$ satisfies the governing equation whereas it remains its initial condition ahead of $x_{\mathrm s}$:  
\begin{eqnarray} \label{eqn:sto_shock}
	\Pi(K, x,t)= \left \{
	\begin{array} {lll}
	 \mathcal H( K - z\sin(\pi tx)),  & \quad x>x_{\mathrm s}(t)\\[10pt]
	\mathcal H(K - k ),  &\quad x<x_{\mathrm s}(t) 
	\end{array} \right..
\end{eqnarray}

Following the proposed CDF scheme, we calculate the cumulative density function $F_{\mathrm k}$ from solutions of $\Pi$. Figure~\ref{fig:BGsolution} presents the $l_2$-norm errors~\eqref{error:1D-random} between the converged solution $F_{\mathrm k}(K; x=0.4, \, t=1)$ and those from CDF method, and the MCS simulations of the Burger's equation, respectively, at different realization numbers $M$. Here the reference solution is obtained from $2000$ MCS simulations. 

\begin{figure}[H]
\centering
\includegraphics[height=2.1in]{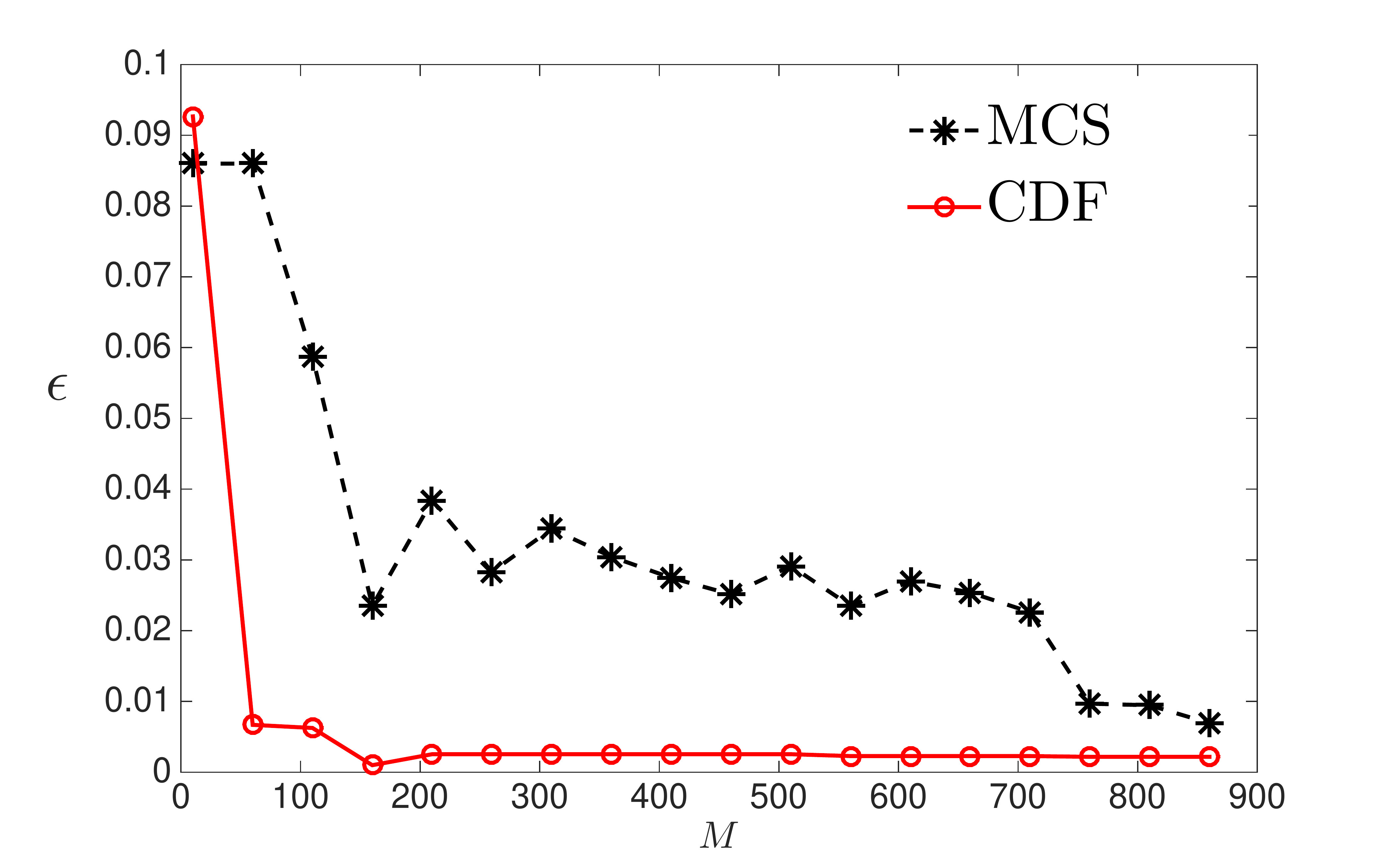}
  \caption{The $l_2$-norm errors~\eqref{error:1D-random} between the converged solution $F_{\mathrm k}(K; x=0.4, \, t=1)$ and those from CDF method (CDF), and those from Monte Carlo simulations (MCS) of the Burger's equation, respectively, at different realization numbers $M$. The converged solution is obtained from $2000$ MCS simulations.}
\label{fig:BGsolution}
\end{figure}
\end{revision2}

%%%%%%%%%%%%%%%%%%%%%%%%%%%%%%%%%%%%%%%%%%%%%%%%%%%%%%%%
%%%%%%%%%%%%%%%%%%%%%%      RESULTS      %%%%%%%%%%%%%%%%%%%
%%%%%%%%%%%%%%%%%%%%%%%%%%%%%%%%%%%%%%%%%%%%%%%%%%%%%%%%
%%%%%%%%%%%%%%%%%%%%%%%%%
\section{Results and discussion} \label{sec:results}
%%%%%%%%%%%%%%%%%%%%%%%%%
In this section, we investigate the robustness of our numerical scheme through an environmental application of the kinematic wave model. To be specific, we consider a one-dimensional Manning open-channel flow under three source cases, namely: no source ($S=0$), constant source ($S=\mbox{const}$) and spatially-dependent source ($S=S(x)$).  Solutions from our CDF approach and those from the direct simulations would be compared and analyzed. 

The Saint-Venant equation is the unidirectional form of shallow water equations and often used to describe flood wave propagation. Using the Manning constitutive relationship, it can be written as:
\begin{equation} \label{def:manning}
	\frac{\partial k}{\partial t} + \frac{\partial q}{\partial x} = S,
	\qquad q = \frac{ \sqrt{s_0} }{C_{\mathrm M}} k^{4/3}, 
\end{equation}
%
%in a flow domain (e.g., the length of a river downstream from $x=0$) of $L = 20\,\mathrm{km}$. 
where $k(x,t)$ $[\mathrm m^2]$ is the cross-sectional area of a channel occupied by the fluid at a point $x$ along the channel length, $q(x,t)$ $[ {\mathrm m}^3 / {\mathrm s}]$ describes the volumetric flow rate, $S$ $[\mathrm m^2 /  \mathrm s]$ is the lateral inflow rate of tributaries and/or upstream rainfall rate, $s_0(x)$ is the slope of the channel bed, and $C_{\mathrm M}(x)$ [s/m$^{1/3}$] denotes the Manning's roughness coefficient. We consider $C_{\mathrm{M}}$ and $s_0$ as random fields. Thus, we have the flux function formulation \eqref{eqn:con}, with $\boldsymbol{z} = (z_1, z_2) = (C_{\mathrm{M}}, s_0)$.% these variable represent the random function $z_{1}(x),z_{2}(x)$.

The Saint-Venant model~\eqref{kwt-saint-venant} provides a good approximation of flood waves for a Froude number smaller than one, in which the main disturbance is carried downstream only by the kinematic waves while dynamic waves (long gravity waves) attenuate rapidly~\cite{Lighthill_1955a_Kinematic}. We also note that the kinematic wave model neglects influence on the river upstream of the junction~\cite{Lighthill_1955a_Kinematic} and the backwater effects (upstream propagation caused by local acceleration, convective acceleration and pressure), the flow rate throughout the flow domain is therefore non-negative, $q(x,t) \geq 0$. Without loss of generality, the approach presented below can also incorporate higher physical dimensions and other types of constitutive relationship, such as the Ch{\' e}zy formula, to represent a balance between the friction at the channel bottom and the gravitational force.

With no analytical solution available, the Saint-Venant equation~\eqref{kwt-saint-venant} is often solved numerically in the form of flux $q$~\cite{wang-2012-uncertainty}:% in a conservation form:
\begin{equation} \label{kwt-saint-venant}
	\frac{3}{4} \, \left(\frac{C_{\mathrm M}}{\sqrt{s_0}} \right) ^{\frac{3}{4}} q^{-\frac{1}{4}} \, \frac{\partial q}{\partial t} + \frac{\partial q}{\partial x} = S,
\end{equation}
Such equation would be directly solved via a fifth-order weighted essentially non-oscillitorary (WENO) scheme in space and third-order TVD Runge-Kutta scheme (Appendix~\ref{appendix:weno}). 

In order to compare with the direct simulation of the hyperbolic equation above~\eqref{kwt-saint-venant}, we apply the CDF formulation to solve for flux $q$. The corresponding fine-grained CDF equation of flux, $\Pi = \mathcal H\left[ Q - q(x,t)\right]$ is:
\begin{equation} \label{eqn:pi-saint-venant}
	\frac{3}{4} \, \left(\frac{C_{\mathrm M}}{\sqrt{s_0}} \right) ^{\frac{3}{4}} Q^{-\frac{1}{4}} \, \frac{\partial \Pi}{\partial t} + \frac{\partial \Pi}{\partial x}  + S \frac{\partial \Pi}{\partial Q} =0.
\end{equation}
%

%Substituting the Manning formula~\eqref{def:manning}, it can be written in a conservation form:
%% 
%\begin{equation} \label{kwt-saint-venant}
%	\frac{\partial k}{\partial t} + \frac{4 \sqrt{s_0} k ^{1/3}}{C_{\mathrm M}}  \, \frac{\partial x}{\partial t} + \frac{k^{4/3}}{2C_{\mathrm M} \sqrt{s_0}} \frac{\partial s_0}{\partial x} -\frac{k^{4/3}\sqrt{s_0}}{C_{\mathrm M}^2 } \frac{\partial C_{\mathrm M}}{\partial x} = S,
%\end{equation}
%%

Following statistical data from earlier analysis~\cite{Buhman-2002-Stochastic,Gates_1996a_Spatio,Gates_1996b_Spatio,Moramarco_2000_Kinematic}, both slope $s_0(x)$ and the Manning coefficient $C_{\mathrm M}(x)$ are treated as stationary random fields exhibiting lognormal distribution and exponential correlation function with correlation length $\lambda$. The mean and standard deviation are $0.01$, $0.0025$ and $0.037$ [s/m$^{1/3}$], $0.00925$ [s/m$^{1/3}$], for $s_0(x)$ and $C_{\mathrm M}(x)$, respectively. We employ Gaussian process model to produce $M=1000$ realizations for each random parameters. Flow rate at initial time and the inlet is set as $q_{\mathrm{in}}(x, t=0)=0.5$ [m$^3$/s] and $q_{\mathrm{bx}}(x=0, t) = \mathrm{max}\left( \sin \pi t, 0.5\right)$ [m$^3$/s], respectively.

\begin{figure}[hbtp]
\centering
\includegraphics[height=2in]{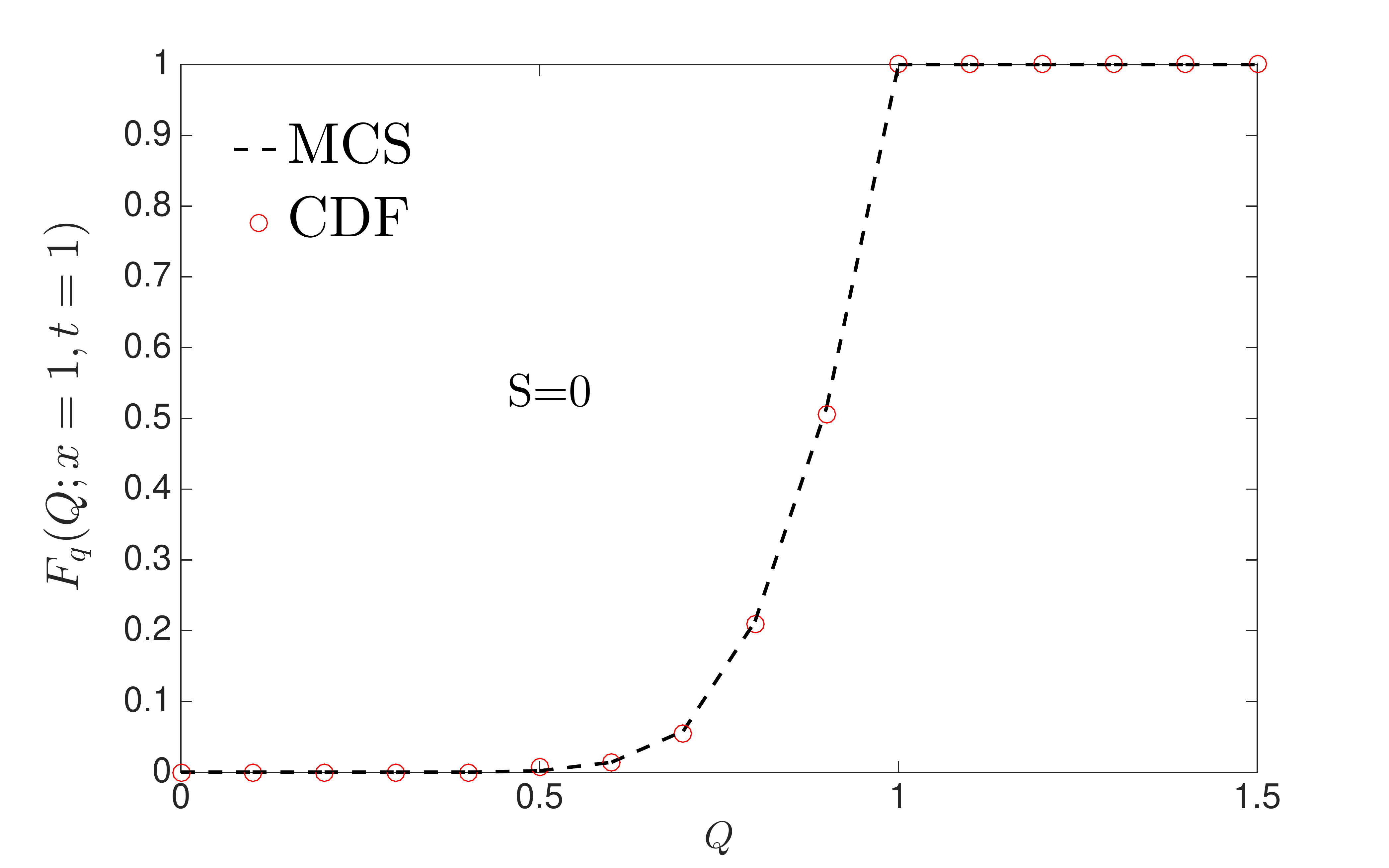}
\includegraphics[height=2in]{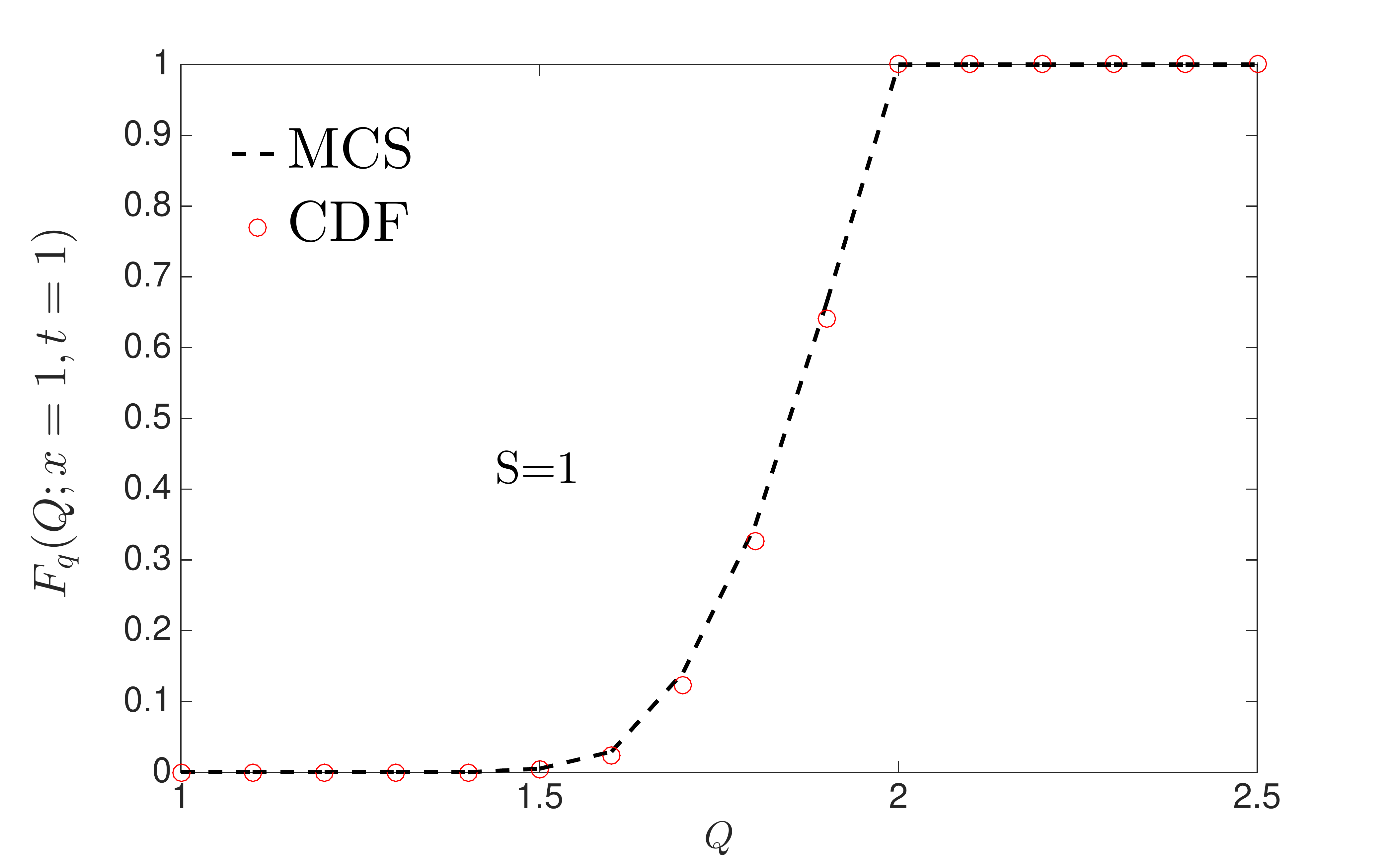}
\includegraphics[height=2in]{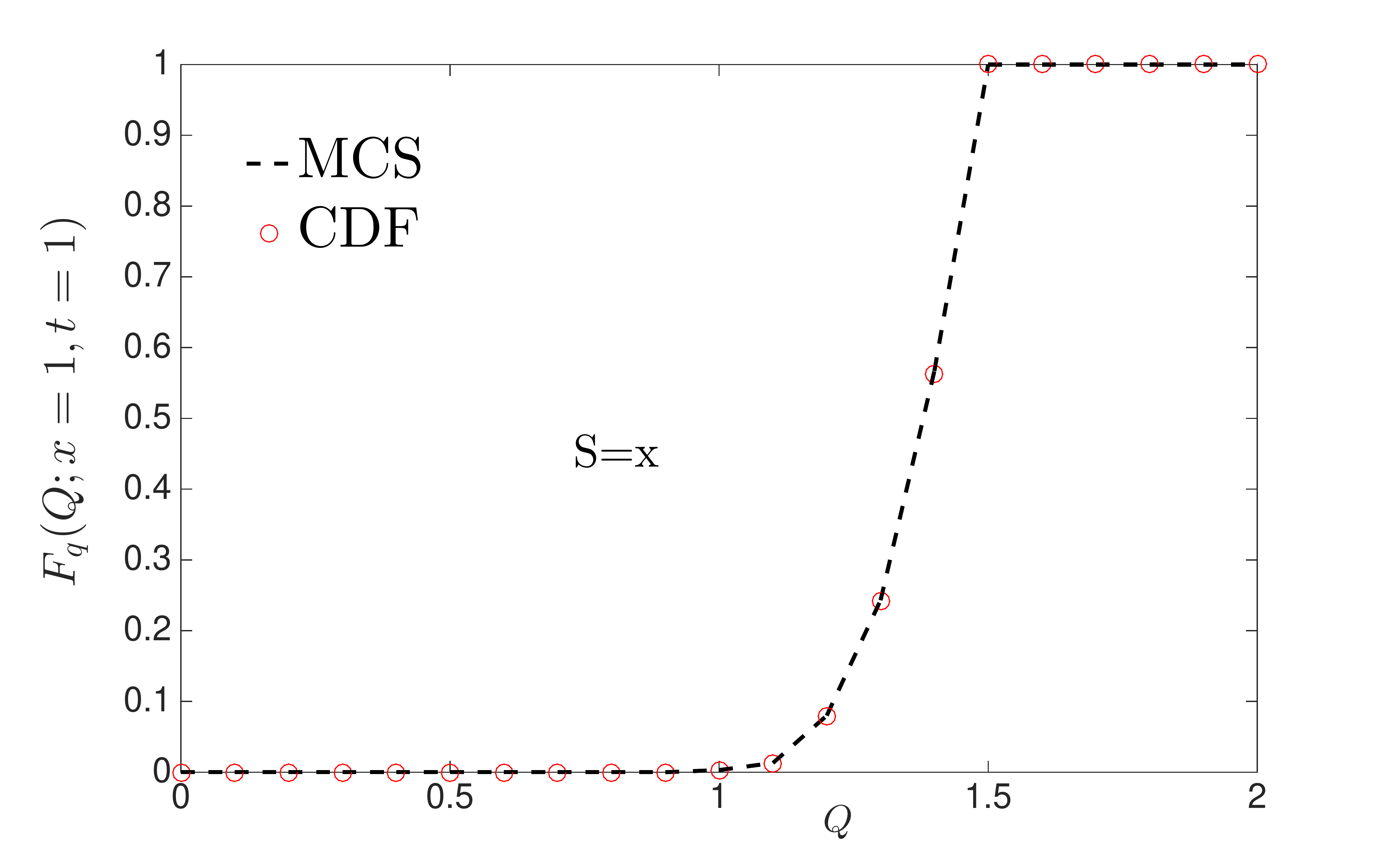}
  \caption{Cumulative density function of flux $F_{\mathrm q}(Q; x=1, \, t=1)$ from CDF scheme (CDF) and Monte Carlo Simulations (MCS) of the kinematic wave equation~\eqref{kwt-saint-venant}, under $S=0$, $S=1$ and $S=x$ conditions.}
\label{fig:s-solution}
\end{figure}

We conduct $1000$ MCS simulations for the original kinematic wave equation~\eqref{kwt-saint-venant} and use it as the reference solutions. Figure~\ref{fig:s-solution} illustrates the cumulative density function of flux $F_{\mathrm q}(Q; x=1, \, t=1)$, under three source conditions: $S=0$, $S=1$ and $S=x$. It is clear that solutions from the CDF method match well with those from the direct MCS simulations. For such desired ensemble accuracy, the fine-grained CDF equation at each realization requires significantly less CPU time to compute (Table~\ref{tab1:speed}). This is one of the main advantages of our approach using characteristics: although using direct solvers for the kinematic wave equation (in this case WENO-Roe finite-difference schemes) seems straightforward, using the method of characteristics by first deriving the $\Pi$ equation results in more than 50X speedup with a negligible compromise in accuracy.

\begin{table}[H]
\centering
\begin{tabular}{lccc}
\hline
                 &   $S=0$   &$S=1$     &$S=x$  \\
\hline
MCS (seconds)      &46.22   &45.05   & 45.99        \\
CDF (seconds)     &2.00  &1.98     &2.00        \\
\hline
\end{tabular}
\caption{The CPU time of the direct simulation (MCS) and the fine-grained CDF approach (CDF) for a single realization, using a workstation with an Intel (R) Core (TM) i5-4210H 2.90 GHz with 8.00 GB RAM.}
\label{tab1:speed}
\end{table}

\begin{figure}[hbtp]
\centering
\includegraphics[height=2in]{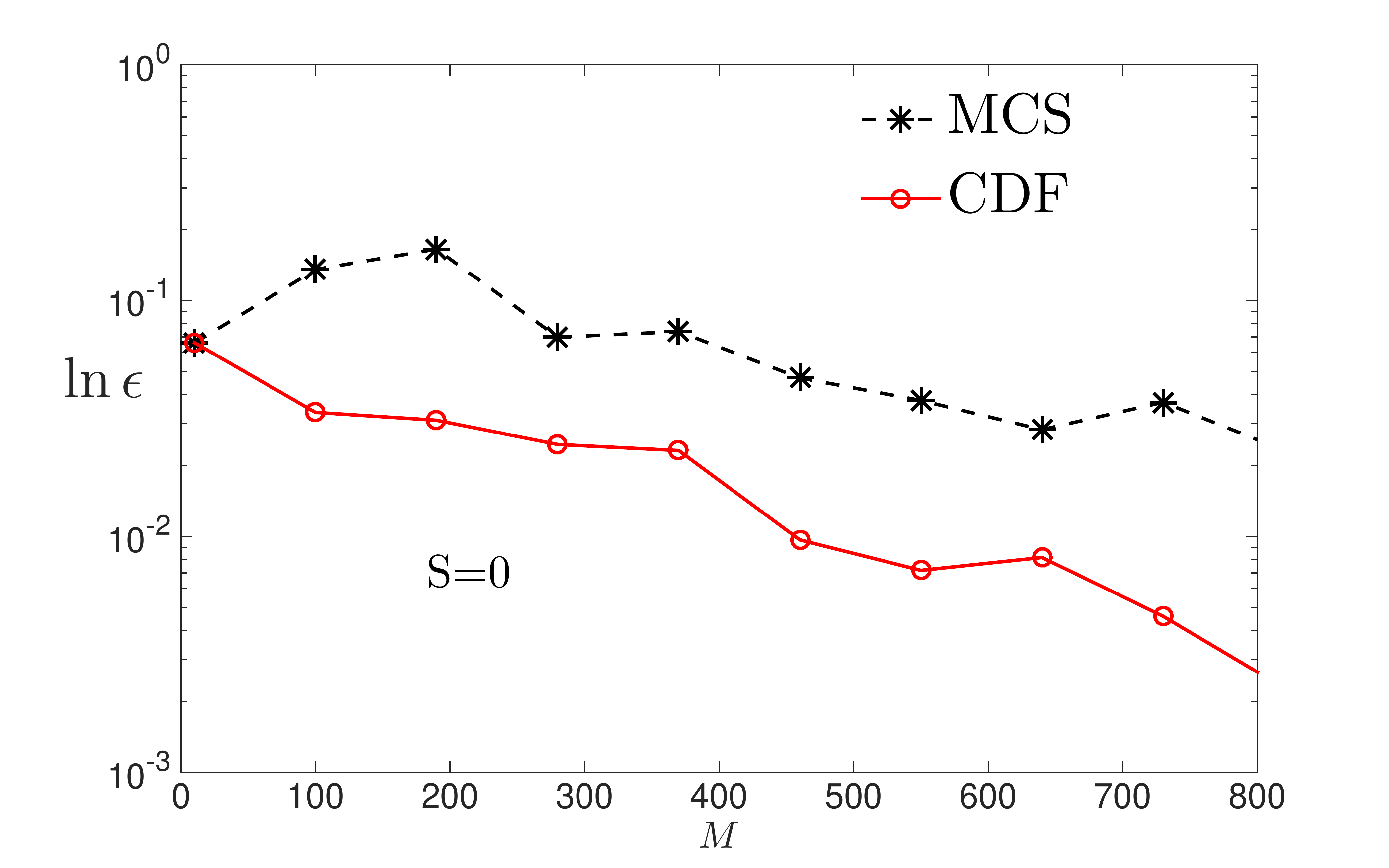}
\includegraphics[height=2in]{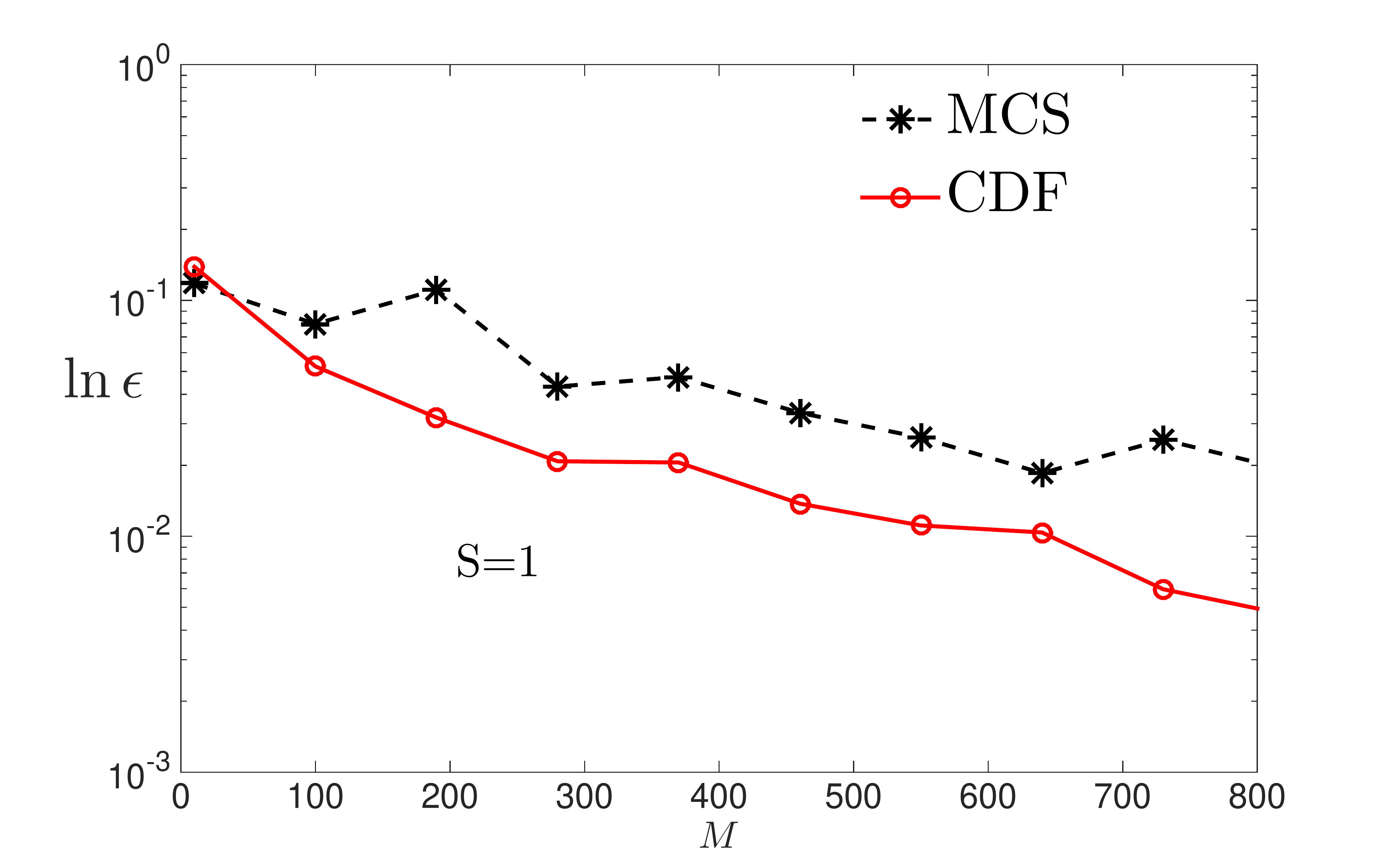}
\includegraphics[height=2in]{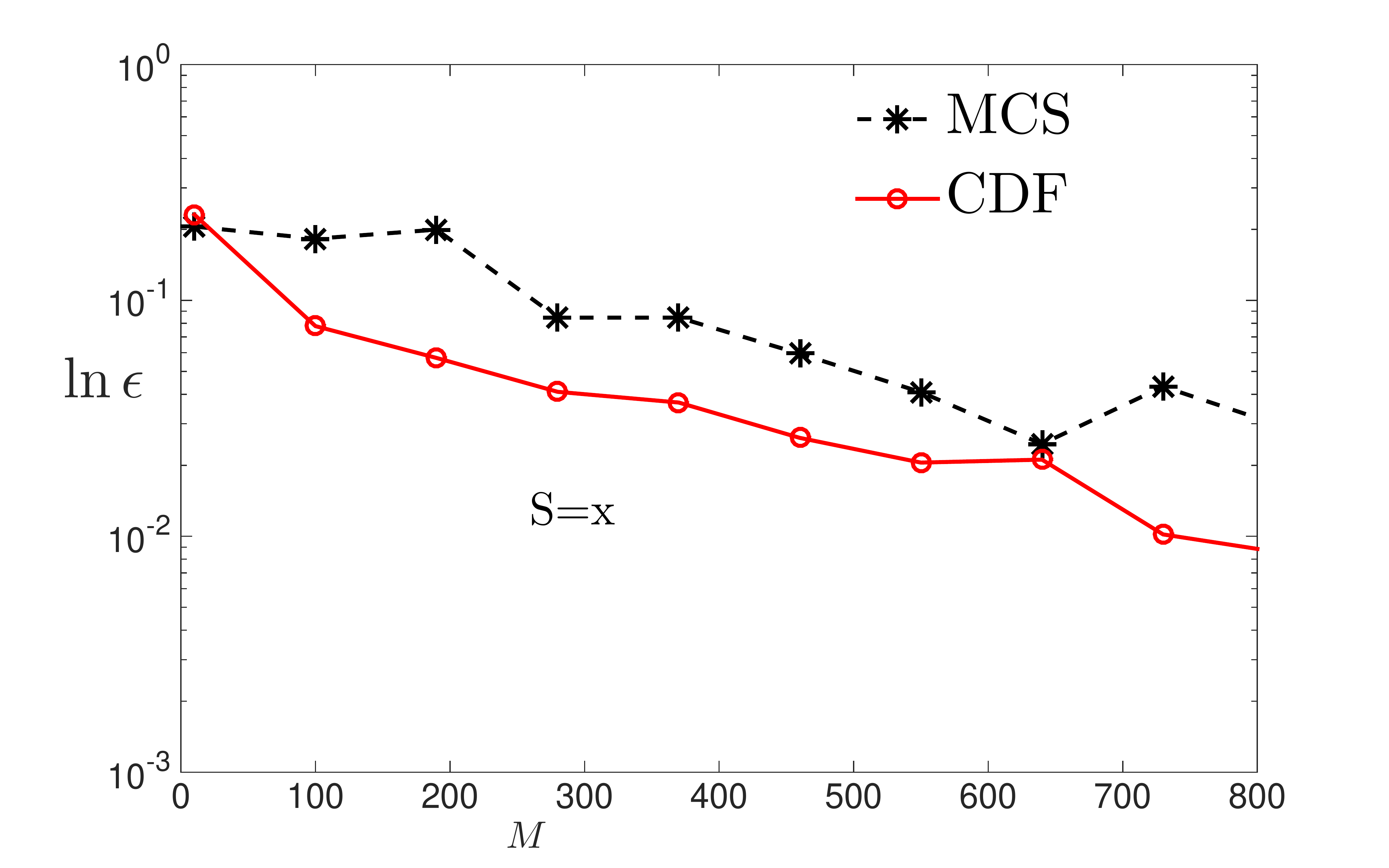}
  \caption{The relative errors $\ln \epsilon$~\eqref{error:1D-random} between benchmark solution $F_{\mathrm q}(Q; x=1, \, t=1)$ and those from the CDF scheme, and from MCS, respectively, at different realization numbers $M$, under $S=0$, $S=1$ and $S=x$ conditions. The benchmark solution here is taken as the MCS solution from $1000$ realizations.}
\label{fig:s-solution-convergence}
\end{figure}

We now examine the convergence rates of the CDF method and direct MCS. Their relative error~\eqref{error:1D-random} to the benchmark MCS solution $F_{\mathrm q}(Q; x=1, \, t=1)$ are plotted in Fig.~\ref{fig:s-solution-convergence} for different realization numbers $M$. We find that the proposed CDF scheme yields smaller errors than those from the MCS for the same number of realizations and overall exhibits faster convergence rate.
%A relative error is introduced: 
%%
%\begin{equation} \label{def:error}
%	|| \epsilon ||_{2}=\frac{\|F_{\mathrm q}^{(i)}- F_{\mathrm{ref}}\|_{2}}{ \|F_{\mathrm{ref}}\|_{2} }, \quad i \leq 1000,
%\end{equation}
%%
%in which $F_{\mathrm q}^{i}$ denotes the ensemble result of $i$ realizations using the CDF method (or the direct simulations), $F_{\mathrm{ref}}$ is the CDF solution from $1000$ realizations using the same numerical method; and $\| \cdot \|_{2}$ refers to the second norm. 

%
%

%\begin{comment}
%Lastly, we investigate the degree of parametric uncertainty on the convergence rate of our numerical scheme. Introducing $\gamma(x) = (\sqrt{s_0} / C_{\mathrm M} ) ^{-3/4}$ as a joint representation of the random parameters, $s_0$ and $C_\mathrm{M}$, we describe the degree of parametric uncertainty through the coefficient of variation $CV_{\gamma}$ (absolute value of the ratio of the standard deviation to the mean). In Fig.~\ref{fig:variance}, the relative error $\epsilon$ of $F_{\mathrm q}(Q; x=1, \, t=1)$ are presented at various realization number $M$. One can see that as the degree of uncertainty rises from $CV_{\gamma}=0.01$ to $0.25$, convergence rate drops in tandem and more realizations are needed. In other words, at relative deterministic case ($CV_{\gamma}=0.01$), a small number of realizations (around three hundred) are sufficient to obtain the statistical distribution of system state.
%
%%
%\begin{figure}[h]
%\centering
%%\includegraphics[width=3in]{figures/CV-001.eps}
%%\includegraphics[width=3in]{figures/s0-convergence.eps}
%\includegraphics[height=2in]{figures/CV-compare}
%  \caption{Relative error $\epsilon$ of $F_{\mathrm q}(Q; x=1, \, t=1)$ at different realization number for: $CV_{\gamma}=0.01$ and $0.25$, respectively. Both are computed using the CDF scheme.}
%\label{fig:variance}
%\end{figure}
%%
%\end{comment}
%%%%%%%%%%%%%%%%%%%%%%%%%%%%%%%%%%%%%%%%%%%%%%%%%%%%%%%%
%%%%%%%%%%%%%%%%%%%%%%      CONCLUSION       %%%%%%%%%%%%%%%%%%%
%%%%%%%%%%%%%%%%%%%%%%%%%%%%%%%%%%%%%%%%%%%%%%%%%%%%%%%%

%%%%%%%%%%%%%%%%%%%%%%%%%%%%%%%%%%%%%%%%%%%%%%%%%%%%%%%%
\section{Conclusion}\label{sec:summary}
%%%%%%%%%%%%%%%%%%%%%%%%%%%%%%%%%%%%%%%%%%%%%%%%%%%%%%%%

In this paper, we presented a numerical framework to implement the CDF method in order to obtain the cumulative density function (CDF) of the system state in the kinematic wave model. The approach relies on solving the fine-grained CDF equation derived from the CDF method. Its accuracy and robustness were investigated via comparison with direct MCS of several numerical examples and a kinematic wave system, the one-dimensional Saint-Venant equation. Our analysis leads to the following conclusions: 

\begin{itemize}

\item{In contrast to previous CDF approach that focuses on the derivation and computation of an ensemble-averaged equation of system state distribution, we directly solve the fine-grained CDF equation, which is exact from the original stochastic system.}

\item{At each single realization, the proposed numerical scheme proves to be computationally more efficient in solving the linear fine-grained CDF equation than the   direct simulation of the nonlinear kinematic wave equation.}

\item{In obtaining the cumulative density function of the system state, our scheme exhibits superior convergence rate and requires fewer realizations than the direct MCS.}

%\item{More variability in input uncertainty would require more realizations for the CDF approach to converge with acceptable accuracy. This is a limitation of the Monte Carlo estimator we use for the CDF, and not a limitation introduced by our approach.}

\item{Faster convergence rate can be achieved with improved quadrature rules and would be the focus of future works.} 

%\item{Instead of approximating the CDF integral using quadrature rules, one may also employ MCS of the fine-grained CDF equation and thus tackle systems with random parameters described as random fields and thus does not impose a prior on the number, distribution type or covariances of random parameters .}

\end{itemize}
%%%%%%%%%%%%%%%%%%%%%%%%%%%%%%%%%%%%%%%%%%%%%%%%%%%%%%%%
%%%%%%%%%%%%%%%%%%%%%%      APPENDIX       %%%%%%%%%%%%%%%%%%%
%%%%%%%%%%%%%%%%%%%%%%%%%%%%%%%%%%%%%%%%%%%%%%%%%%%%%%%%
\appendix

%%%%%%%%%%%%%%%%%%%%%%%%%%%%%%%%%%%%%%%
\section{Derivation of Raw CDF Equation} \label{appendix:pi}
%%%%%%%%%%%%%%%%%%%%%%%%%%%%%%%%%%%%%%%

Following the definition of $\Pi$ in \eqref{def:Pi}, we find its spatial and temporal derivatives as
\begin{eqnarray} \label{eqn:pixt}
	\frac{\partial \Pi}{\partial t} = \frac{\partial \Pi}{\partial k}\, \frac{\partial k}{\partial x} =-\frac{ \partial \Pi}{\partial K} \frac{\partial k}{\partial t},
	\qquad  \nabla \Pi = \frac{\partial \Pi}{\partial k} \nabla k = - \frac{ \partial \Pi}{\partial K} \nabla k,
\end{eqnarray}

For smooth solutions, we multiply the kinematic wave equation~\eqref{KWT} with $\partial \Pi / \partial K$ and substitute the derivates above~\eqref{eqn:pixt}:
\begin{eqnarray}\label{eq:A2}
	\frac{ \partial \Pi } {\partial t} +  
	 \frac{ \partial \Pi } {\partial K} \left( \frac{\partial \mathbf q}{\partial k}  \cdot \nabla k +
	  \nabla_{\mathbf z} \mathbf q \cdot \nabla \mathbf{z} -S \right) =0,
\end{eqnarray}
where the operator $\nabla_{\mathbf z}(\cdot)$ is defined as:
\begin{eqnarray}
	 \nabla_{\mathbf z} \equiv  (\frac{\partial}{\partial z_1}, \,..., \, \frac{\partial}{\partial z_{\mathrm N} } )^\mathrm{T},
\end{eqnarray}

We note here that derivative of the fine-grained CDF in the probability space is the dirac delta function, e.g. $\partial \Pi / \partial K = \delta (K-k)$. Using its sifting property, $g(k) \delta(K-k)= g(K) \delta(K-k)$ for any test function $g(\cdot)$, all $k(\mathbf x, t)$ in the equation~\eqref{eq:A2} can be replaced by $K$ and leads to a linear fine-grained CDF equation~\eqref{eqn:pi}. Interested readers can refer to earlier studies~\cite{Pope2000_turbulent, wang-2012-uncertainty} for detailed derivations.

%%%%%%%%%%%%%%%%%%%%%%%%%%%%%%%%%%%%%%%
\section{TVD Runge-Kutta scheme for characteristic lines}
\label{appendix:RK}
%%%%%%%%%%%%%%%%%%%%%%%%%%%%%%%%%%%%%%%
The Runge-Kutta method is a simple and robust iterative scheme. Here we employ a third-order TVP Runge-Kutta scheme~\cite{Gottlieb-1998-Total} to solve the ordinary differential equations~\eqref{eqn:cl-v} and solution at the $h$-th stencil and the $n$-th time-step, $k_{h}^n $, can be obtained in a three-step iteration: 
\begin{eqnarray*}
  \mathbf c^{n1}_{h} &=& \mathbf c^{n-1}_{h}+\Delta t \, \mathbf v \left[ \mathbf c^{n-1}_{h},S(t^{n}) \right], \\
  \mathbf c^{n2}_{h} &=& \frac{3}{4} \mathbf c^{n-1}_{h}+\frac{1}{4} \mathbf c^{n1}_{h}+\frac{1}{4} \, \Delta t \,  \mathbf v \left[ \mathbf c^{n-1}_{h},S(t^{n}+\Delta t) \right], \\
  \mathbf c^{n}_{h} &=& \frac{1}{3} \mathbf c^{n-1}_{h}+\frac{2}{3} \mathbf c^{n2}_{h}+\frac{2}{3} \, \Delta t \,  \mathbf v \left[ \mathbf c^{n-1}_{h},S(t^{n}+\frac{1}{2}\Delta t) \right], 
\end{eqnarray*}
%

%%%%%%%%%%%%%%%%%%%%%%%%%%%%%%%%%%%%%%%
\section{WENO-Roe Numerical scheme for kinematic wave equation}
\label{appendix:weno}
%%%%%%%%%%%%%%%%%%%%%%%%%%%%%%%%%%%%%%%
The WENO scheme~\cite{Shu-1998-Essentially} provides a nonlinear adaptive procedure to automatically select the smoothest local stencil in numerical approximation of fluxes. We first take the change of variable, $q' = ( q \,C_\mathrm{M}  / \sqrt{s_0})^{3/4}$, and rewrite the one-dimensional kinematic wave equation~\eqref{kwt-saint-venant} as:
\begin{eqnarray} \label{eqn:weno}
	\frac{\partial q'}{\partial t} + \frac{\partial }{\partial x} \left[ \left( \frac{\sqrt{s_0}}{C_\mathrm{M}} \right) q'^{\frac{4}{3} }\right] = S,
\end{eqnarray}

The spatial domain $x \in [a, b]$ is divided into $N_\mathrm{x}$ stencils of size $\Delta x$ and let us denote $q'_h$ as the value of $q'$ at the $h$-th stencil. Using a fifth-order finite difference WENO-Roe scheme, we can numerically approximate the partial differential equation with an ordinary differential equation:
\begin{subequations}
\begin{eqnarray} \label{eqn:weno-x}
	&& \qquad 
	\frac{\mathrm d q'_h}{\mathrm d t} = - \frac{1}{\Delta x} \left[ \hat{u}(q')_{h+\frac{1}{2}} - \hat{u}(q')_{h-\frac{1}{2}} \right] + S_h, \\
	&& \qquad
	\hat{u}(q')_{h+\frac{1}{2}} = \begin{cases}
	u(q')^-_{h+\frac{1}{2}}, \quad \mbox{if}\, a_{h+\frac{1}{2}} > 0\\
	u(q')^+_{h+\frac{1}{2}}, \quad\mbox{if} \, a_{h+\frac{1}{2}} < 0
	\end{cases}, \\
	&&	\qquad
	 \hat{u}(q')_{h -\frac{1}{2}} = \begin{cases}
	u(q')^-_{h-\frac{1}{2}}, \quad \mbox{if} \, a_{h-\frac{1}{2}} > 0\\
	u(q')^+_{h+\frac{1}{2}}, \quad \mbox{if} \, a_{h-\frac{1}{2}} < 0
	\end{cases},
\end{eqnarray}
whose Roe speeds at the $h$-th stencil is:
\begin{eqnarray} \label{def:roe}
	 a_{h+\frac{1}{2}}\equiv\frac{u(q')_{h+1} - u(q')_{h}} {q'_{h+1}-q'_{h}},
	\qquad \qquad
	a_{h-\frac{1}{2}}\equiv\frac{u(q')_{h} - u(q')_{h-1}}{q'_{h}-q'_{h-1}}, 
\end{eqnarray}
\end{subequations}
and $u(q') = \sqrt{s_0}  q'^{4/3} /C_\mathrm{M}$. The superscripts $^{\pm}$ refers to the direction from which one interpolates the half stencils, in other words, if $a_{h+\frac{1}{2}} > 0$, the value $u(q')^-_{h-\frac{1}{2}}$ is approximated from the upwind direction. Exact expressions of those interpolations can be found in earlier studies~\cite{Shu-1998-Essentially}.

Now we can solve the ordinary differential equation with a third-order TVD Runge-Kutta scheme (Appendix~\ref{appendix:RK}). Boundary conditions at the $h_{-2}$-th and $h_{-1}$-th stencils are taken as the same value as that at $h_0$-th.

We validated our WENO-Roe scheme via the deterministic example in section~\ref{subsec:example} with a time step $\Delta t = 0.0001$. Table~\ref{tab1:test-1-weno} shows the comparison between our numerical result and the analytical solution~\eqref{sol:test1} at time $t=0.1$: as the space step $\Delta x$ drops (from $0.05$ to $0.00625$), the 2-norm difference ($\epsilon$) between the two solutions is reduced to $1.97 \times 10^{-9}$ while the WENO-Roe convergence rate is approximately $4.5$, providing reasonable validation of the fifth-order approach.
\begin{table}
\centering
\begin{tabular}{lcc}
\hline
$\Delta x$  &$\epsilon$ & Convergence Rate \\
\hline
0.05          & $2.22 \times 10^{-5} $     & None    \\
0.025         &$1.02 \times 10^{-6} $     & 4.44    \\
0.0125        &$4.59 \times 10^{-8}$      & 4.47    \\
0.00625        &$1.97 \times 10^{-9}$      & 4.54    \\
\hline
\end{tabular}
\caption{Comparison between the WENO-Roe numerical result and the exact solution~\eqref{sol:test1} at time $t=0.1$.}
\label{tab1:test-1-weno}
\end{table}
%

%%%%%%%%%%%%%%%%%%%%%%%%%%%%%%%%%%%%%%%%%%%%%%%%%%%%%%%%%%%%%%%%%%%%%%%%%%%%%%%%
%%%%%%%%%%%%%%%%%%%%%%%%%%%%%%%%%%%%%%%% ACKNOWLEDGMENT %%%%%%%%%%%%%%%%%%%%%%%%%%%%%%%%%%%%%%%%
%%%%%%%%%%%%%%%%%%%%%%%%%%%%%%%%%%%%%%%%%%%%%%%%%%%%%%%%%%%%%%%%%%%%%%%%%%%%%%%%
\section*{Acknowledgment}
P. Wang and M. Cheng were partially funded by the National Natural Science Foundation of China (Grant No. 11571028), the National Key Research and Development Program of China (Grant No. 2017YFB0701702) and the Recruitment Program of Global Experts. A.~Narayan was partially supported by AFOSR FA9550-15-1-0467. X.~Zhu was partially supported by Simons Foundation.

\section*{References}

\bibliographystyle{siam}
\bibliography{Bib/pdf-cdf,Bib/k-wave,Bib/collocation,Bib/approximation,Bib/random}

\end{document}